\documentclass{amsart}

\usepackage[colorlinks,driverfallback=dvipdfm, urlcolor=cyan, linkcolor=blue,unicode,psdextra]{hyperref}
\usepackage{color,graphicx,shortvrb}
\usepackage[utf8]{inputenc}

\usepackage[active]{srcltx} 

\usepackage{tikz-cd}

\usepackage{parskip}

\usepackage{faktor}

\usepackage{ bbm }

\usepackage{comment}

\usepackage{amsmath,amssymb, amsfonts, amsthm, verbatim,amscd,bm,color,yhmath,bm, enumerate}
\usepackage[normalem]{ulem}
\usepackage[driverfallback=dvipdfm]{hyperref}

\usepackage[backend=biber, style=nature, sorting=nyvt,
doi=false,isbn=false,url=false,eprint=false]{biblatex}
\AtEveryBibitem{\clearlist{language}}

\usepackage[active]{srcltx} 

\newtheorem{lem}{Lemma}[section]
\newtheorem{thm}[lem]{Theorem}
\newtheorem{prop}[lem]{Proposition}
\newtheorem{cor}[lem]{Corollary}

\theoremstyle{definition}
\newtheorem{remark}[lem]{Remark}

\newtheorem{problem}[lem]{Problem}

\numberwithin{equation}{section}

\newcommand{\A}{\mathfrak{A}}
\newcommand{\At}{\mathcal{A}}
\newcommand{\Nat}{\mathcal{N}}

\usepackage{ulem}
\usepackage{marginnote}
\usepackage{geometry}

\begin{document}


\title
{Estimations of the numerical index of a JB$^*$-triple}

\author[D. Cabezas]{David Cabezas}
\address[D. Cabezas]{
Departamento de An{\'a}lisis Matem{\'a}tico, Facultad de
Ciencias, Universidad de Granada, 18071 Granada, Spain.}
\email{dcabezas@ugr.es}

\author[A.M. Peralta]{Antonio M. Peralta}
\address[A.M. Peralta]{Instituto de Matem{\'a}ticas de la Universidad de Granada (IMAG), Departamento de An{\'a}lisis Matem{\'a}tico, Facultad de
Ciencias, Universidad de Granada, 18071 Granada, Spain.}
\email{aperalta@ugr.es}


\begin{abstract} We prove that every commutative JB$^*$-triple has numerical index one. We also revisit the notion of commutativity in JB$^*$-triples to show that a JBW$^*$-triple $M$ has numerical index one precisely when it is commutative, while $e^{-1}\leq n(M) \leq 2^{-1}$ otherwise. Consequently, a JB$^*$-triple $E$ is commutative if and only if $n(E^*) =1$ (equivalently, $n(E^{**}) =1$). In the general setting we prove that the numerical index of each JB$^*$-triple $E$ admitting a non-commutative element also satisfies $e^{-1}\leq n(M) \leq 2^{-1}$, and the same holds when the bidual of $E$ contains a Cartan factor of rank $\geq 2$ in its atomic part.
\end{abstract}

\maketitle

\section{Introduction}

This paper is aimed to present the first estimations of the numerical index of the complex Banach spaces associated with JB$^*$-triples and its relation with the commutativity of these objects (see below for definitions). Despite the intense activity developed since the sixties of the previous century to compute the numerical index of certain concrete classes of Banach spaces (in most of cases complex), our current knowledge is practically nonexistent in the case of this deeply studied class of complex Banach spaces.  

Let us begin by recalling the basic notions. If $X$ is a normed space, we shall denote by $\mathcal{B}_X$ and $S_X$ the unit ball and the unit sphere of $X$, respectively. The symbol $B(X,Y)$ will denote the space of bounded linear operators from $X$ into another normed space $Y$. We shall write $B(X)$ for $B(X,X).$ Suppose $X$ is a real or complex normed space and $T$ lies in $B(X)$. The \emph{numerical range} and the \emph{numerical radius} of $T$ are defined by $$V(T):=\left\{ \phi (Tx) :x\in S_X, \phi\in S_{X^*}, \phi (x)=1\right\}\subset \mathbb{C}, \hbox{ and } v(T):=\sup\left\{|\lambda| : \lambda\in V(T)\right\}\in \mathbb{R}_0^+,$$ respectively. Clearly, $v(T)\leq\|T\|$ for every $T\in B(X)$. The numerical radius is a seminorm on $B(X)$. There are examples of real Hilbert spaces $H$ and $T\in B(H)$ with $v(T) =0.$ The \emph{numerical index} of the space $X,$ defined by
$$n(X):=\max\{k\geq 0:k\|T\|\leq v(T),\ \forall T\in B(X)\}=\inf\{v(T):T\in B(X), \|T\|=1\},$$
 is positive precisely when the numerical radius is a norm (equivalent to the operator norm). Pioneering studies on numerical index trace back to Bauer and Lumer in 1960's. A complete historical outline along with the basic results on numerical index can be found in the celebrated monographs by Bonsall and Duncan \cite{BonDunBookI, BonDunBookII}. For the purposes of this note we recollect some of the milestone results. For example, for each $T\in B(X)$, it is known that $v_{B(X)}(T)=v_{B(X^*)}(T^*)$ (see \cite[\S 9, Corollary 6]{BonDunBookI}). It follows from this that $n(X^*)\leq n(X)$. The values of the numerical indexes of all real Banach spaces exhaust all possibilities in the interval $[0,1]$, while for complex Banach spaces the values cover all elements in $[e^{-1}, 1]$ (cf. \cite{DuncanMcGregorPryceWhite1970}; we note that the below bound in the complex setting is essentially due to the Bohnenblust--Karlin theorem). 

It was an intriguing question, open during decades, whether the inequality $n(X)\geq n(X^*)$ can be strict or not. A definitive answer was given by Boyko, Kadets, Martín, and Werner in \cite{BoykoKadMar2007numIndexDual} who found an example of a Banach space whose numerical index is strictly greater
than the numerical index of its dual (see also \cite[Example 4.3]{MarJFA2008}). The references \cite{MartinSurveyNumIndex2000, KadMarPay2006surveyRACSAM} contain more recent and up-to-date surveys on the numerical index. There is a very intense recent activity led towards the calculation of the numerical indices of several prominent Banach spaces. The numerical index of the $c_0$-, $\ell_1$-, and $\ell_{\infty}$-sum of an arbitrary family of Banach spaces coincides with the infimum of the numerical indices of the respective summands, $n(C(K,X)) = n(X)$, $L^1(\mu,X) = n(X),$ and $L^{\infty}(\nu,X) = n(X),$ for any compact Hausdorff space $K$, any positive measure $\mu$, every $\sigma$-finite measure $\nu,$ and any Banach space $X$ (see \cite{MarPay2000Studia, MarVillLinfty2003}). Moreover, if $M$ is a subspace of $C(K, X)$ which satisfies the property $J$ and is left multiplication invariant, then
$n(M) \leq n(X)$ \cite{BakerBotelho2022}. The numerical index of each function algebra is $1$ \cite{WernerJFA1997}. The numerical index of the projective and the injective tensor product of two Banach spaces is less than or equal to the minimum of the numerical indexes of both factors \cite{MarMerQuerIndexTensorProd2021}. Many other recent publications have been devoted to the study of the numerical index (see, for example, \cite{MartProcAMS2009Lembedded, ChicaMarMer2014,ChicMer2015,KadMarMerPerQue2020, MerQueLMA2021}, there are over 1600 papers published under the MSC item 47A12
``Numerical range, numerical radius'').

One of the most renowned, and cited, results is a result by Huruya affirming that the numerical index of  a C$^*$-algebra $A$ is $1$ if $A$ is commutative and $\frac12$ otherwise \cite{HuruyaPAMS1977}. The arguments by Huruya actually rely on a previous contribution by Crabb, Duncan, and McGregor on the algebra numerical index of a C$^*$-algebra \cite{crabb74characterizations}, the Russo--Dye theorem, and a celebrated result by Kaplansky characterizing the commutativity of a C$^*$-algebra by the nonexistence of a non-zero $2$-nilpotent element (cf. \cite[2.12.21]{diximer77C*-algebras}). 

Subsequent studies addressed the question of determining the numerical index of more general non-associative complete normed algebras. Rodríguez-Palacios proved in  \cite[Theorem 26]{RodPal1980} that if $\A$ is a (unital) non-commutative Jordan $V$-algebra, we have $n(\A) = 1$ if $\A$ is associative and commutative, and $n(\A) = \frac12$ otherwise. The same author in collaboration with Iochum and Loupias showed that the same conclusion holds when $\A$ is a (non-necessarily commutative) JB$^*$-algebra \cite{iochum89commutativity}. The reader can take a look at \cite{KadMorRodPal2001}, and \cite[Proposition 3.5.44 and comments in \S 2.1.47 and page 422]{CabRodBookV1} for a more recent approach. 

C$^*$-algebras and JB$^*$-algebras are strict subclasses of the wider family of complex Banach spaces known as JB$^*$-triples. JB$^*$-triples were introduced to classify bounded symmetric domains in complex Banach spaces of arbitrary dimension by Kaup \cite{kaup83riemann} with the aim of extending the classical Riemann mapping theorem. The geometric, algebraic and holomorphic properties of JB$^*$-triples have been intensively studied during the last forty years. However, we completely lack of any estimation of the numerical index of a general JB$^*$-triple. Martín already posed the question whether, as in the case of C$^*$-algebras and von Neumann preduals, the condition $n(Y) = 1$ implies $n(Y^*)=1$ when $Y$ is a JB$^*$-triple or the predual of a JBW$^*$-triple (cf. \cite[comments after Proposition 3.5]{MarMathNach2008}). The problem was affirmatively solved by Martín in the case of JBW$^*$-triple preduals in \cite[Theorem 2.1]{MartProcAMS2009Lembedded}.  

This paper presents the first known estimations on the numerical index of an arbitrary JB*-triple. Surprisingly, despite of  the lacking of a binary product, the key algebraic property to compute the numerical index is ``commutativity''. We recall that a JB$^*$-triple $E$ with triple product $\{.,.,\}$ is commutative if the operators of the form $L(a,b): x\mapsto \{a,b,x\}$ ($a,b\in E$) commute in $B(E)$. We shall revisit the Gelfand theory for commutative JB$^*$-triples in section~\ref{sec: commutativity}. The notion of commutativity in JB$^*$-triples has been considered by consolidated experts like Kaup \cite[\S 1]{kaup83riemann}, Dineen and Timoney \cite{dineen88centroid}, and Friedman and Russo \cite{FriRuCommutative} (the references can be complemented with \cite[\S 4.2.1]{CabRodBookV1} and \cite[\S 3]{BurPeRaRu2010}). We have already recalled Kaplansky's result stating that a C$^*$-algebra is commutative if and only if it has no non-trivial nilpotent elements \cite[2.12.21]{diximer77C*-algebras}. Motivated by this characterization, it was conjectured that a JB$^*$-algebra is associative if and only if it contains no non-trivial nilpotent elements. The conjecture was finally proved to be true by Iochum, Loupias and Rodríguez-Palacios in \cite{iochum89commutativity}. It is not clear how nilpotency can be applied in the setting of JB$^*$-triples (specially because by the extended Gelfand--Naimark axiom $\|\{a,a,a\}\| = \|a\|^3$ for every element $a$ in a JB$^*$-triple $E$). We turn our point of view to the inner ideals generated by a single element which are naturally equipped with a structure of JB$^*$-algebra. In Theorem \ref{thm: commutativity characterizations JB*-triples} we prove some new characterizations of commutativity in JB$^*$-triples, showing that a JB$^*$-triple $E$ is commutative if and only if one of the next statements holds:\begin{enumerate}[$(a)$]
\item For each $a$ in $E$ the image of the operator $L(a,a) = \{a,a,.\}$ is contained in the inner ideal of $E$ generated by $a$, that is, in $E(a) = \overline{\{a,E,a\}}^{\|.\|}$.
\item The atomic part of $E^{**}$ contains no Cartan factors with rank $\geq 2$ nor rank-one type 1 Cartan factors (i.e. Hilbert spaces) of dimension $\geq 2$.
\end{enumerate} The following consequence, obtained in Corollary \ref{cor: nilpotent in non-commutative JB*-triple}, asserts that a JB$^*$-triple $E$ is non-commutative if and only if one of the following statements holds: \begin{enumerate}[$(1)$]
        \item There is an element $a$ in $E$ such that the inner ideal $E(a)$ contains a non-zero $2$-nilpotent element $b$ as JB$^*$-algebra {\rm(}i.e., $\{b,r(a),b\}=0${\rm)}, equivalently, $E(a)$ is non-associative.
        \item Every single generated inner ideal of $E$ is an associative JB$^*$-algebra and the atomic part of $E^{**}$ reduces to a $\ell_{\infty}$-sum of Hilbert spaces and at least one of them has dimension greater than or equal to $2$. 
    \end{enumerate}

Our main conclusions are contained in section~\ref{sec:num index JBstartriples}. We first prove that every commutative JB$^*$-triple has numerical index one (see Lemma~\ref{l the numerical index of a commutative JB*-triple is 1}). Assuming that $M$ is a JBW$^*$-triple (i.e. a JB$^*$-triple which is also a dual Banach space) we show that $n(M)=1$ if $M$ is commutative and $\frac1e \leq n(M) \leq \frac12$ otherwise (cf. Theorem~\ref{t numerical index in JBW*-triples}). This gives a positive solution to the problem posed by Martín in \cite[comments after Proposition 3.5]{MarMathNach2008} in the case that $Y$ is JBW$^*$-triple, since $n(Y)=1 \Leftrightarrow $ $Y$ is commutative $\Leftrightarrow Y^{**}$ is commutative $\Leftrightarrow n(Y^{**}) =1 \Leftrightarrow n(Y^{*}) =1$.  

In order to compute the numerical index of a non-commutative JB$^*$-triple $E,$ we establish that if $E$ contains an element $a$ such that the inner ideal of $E$ generated by $a$ is a non-associative JB$^*$-algebra, the numerical index of $E$ is in the interval $[e^{-1}, 2^{-1}]$ (see Proposition \ref{prop for a JBtriple admitting a non-commutative inner ideal}).

In Theorem \ref{t numerical index non-commutative JB*-triples} we prove that a JB$^*$-triple is commutative if and only if $n(E^*) = n(E^{**}) =1.$ Furthermore, in case that $E$ is non-commutative due to the presence of a Cartan factor of rank greater than or equal to $2$ in the atomic part of its second dual, we have $\frac1e \leq n(E)\leq \frac12$. Just one case of non-commutative JB$^*$-triple remains out from the scope of our results. Concretely, suppose $E$ is a non-commutative JB$^*$-triple satisfying that the atomic part of its second dual is an (infinite) $\ell_{\infty}$-sum of rank-one Cartan factors (i.e., complex Hilbert spaces) and at least one of them has dimension $\geq 2$. We do not know whether $n(E)\leq \frac12$. 

\subsection{Background on JB$^*$-triples}\label{subsec:background} \ 

In this subsection we gather some of the basic definitions and results employed in our arguments. A JB$^*$-triple is a complex Banach space $E$ equipped with a continuous triple product $\{\cdot,\cdot,\cdot\}:E\times E\times E\rightarrow E$, which is linear and symmetric in the outer variables and conjugate linear in the middle one, and satisfies the following conditions:
\begin{enumerate}
    \item The operator $L(a,b)$ on $E$ given by $L(a,b)x=\{a,b,x\}$ satisfies 
    $$L(a,b)L(x,y)-L(x,y)L(a,b)=L(L(a,b)x,y)-L(x,L(b,a)y),$$ for all $a,b,x,y\in E$. \hfill (Jordan identity)
    \item For each $a\in E$, $L(a,a)$ is a hermitian operator with non-negative spectrum.
    \item $\|\{a,a,a\}\|=\|a\|^3$ for every $a\in E$. \hfill (extended Gelfand--Naimark axiom)
\end{enumerate}
The class of JB$^*$-triples strictly contains the class of JB$^*$-algebras (cf. \cite[page 525]{kaup83riemann}), which in turn strictly contains all C$^*$-algebras. More concretely, given any two complex Hilbert spaces $H,K$, every closed subspace of $B(H,K)$ which is closed for the triple product \begin{equation}\label{eq C*-triple product} \{x,y,z\} :=\frac12 (x y^* z+ z y^*x), 
\end{equation} is a JB$^*$-triple. These JB$^*$-triples are known as JC$^*$-triples. In particular, the just presented triple product induces a structure of JB$^*$-triple on each C$^*$-algebra. In the case of JB$^*$-algebras, the triple product is given by \begin{equation}\label{eq JB*-algebra triple product} \{x,y,z\} := (x\circ y^*) \circ z + (z\circ y^*)\circ x - (x\circ z)\circ y^*.
\end{equation} 

It is perhaps worth to recall the notion of JB$^*$-algebra. A \emph{Jordan algebra} is a (non-necessarily associative) algebra $\A$ whose product $\circ$ is commutative and satisfies the \emph{Jordan identity}: $$(a \circ
b)\circ a^2 = a\circ (b \circ a^2).$$ If $\A$ admits a norm $\|. \|$ satisfying $\| a\circ b\| \leq \|a\| \ \|b\|$, $a,b\in \A$, we say that $\A$ is a normed Jordan algebra. If the norm is additionally complete $\A$ is called a \emph{Jordan--Banach
algebra}. Each associative Banach algebra $A$ is a Jordan--Banach with Jordan product $a\circ b =\frac12 ( a b + ba)$. A \emph{JB$^*$-algebra} is a complex Jordan--Banach algebra $\A$ equipped with an algebra involution $^*$ satisfying  $\|U_a ({a^*})\| = \| 2 (a\circ a^*) \circ a - a^2 \circ a^*\|= \|a\|^3$, for all $a\in \A$. For each element $a$ in a Jordan--Banach algebra $\A$ we denote by $M_a$ the Jordan multiplication operator by the element $a$, that is, $M_a(x)=a\circ x$. Classical references on JB$^*$-algebras can be found in the monographs \cite{hanche-olsen84jordan,AlfShultzBook2003,CabRodBookV1}.


A key property of JB$^*$-triples, known as Kaup's Banach--Stone theorem, asserts that a linear bijection between JB$^*$-triples is an isometry if and only if it is a triple isomorphism (cf. \cite[Proposition 5.5]{kaup83riemann}). We deduce from this deep result that the above expressions of triple products are essentially unique. 

Given $a,b\in E$, we shall write $Q(a,b)$ for the conjugate linear operator on $E$ defined by $Q(a,b) = \{a,x,b\}$ ($x\in E$). The symbol $Q(a)$ will stand for $Q(a,a)$.


A JB$^*$-triple $E$ is called a JBW$^*$-triple if it has a predual. A result by Barton and Timoney assures that each JBW$^*$-triple admits a unique isometric predual and its triple product is separately w$^*$-continuous (see \cite[Theorem 1.4]{barton86weak}). The bidual, $E^{**}$, of each JB$^*$-triple $E$ is a JBW$^*$-triple for a certain triple poduct extending the one in $E$ (cf. \cite[Corollary 11]{dineen86complete}. 


A \emph{tripotent} in a JB$^*$-triple $E$ is an element $e\in E$ satisfying $\{e,e,e\}=e$. Each tripotent $e\in E$ produces a \emph{Peirce decomposition} of $E$ in the form $E=E_2(e)\oplus E_1(e)\oplus E_0(e)$, where $E_k(e):=\{x\in E:L(e,e)x=\frac{k}{2}x\}$ is called the Peirce-$k$ subspace. For $k=0,1,2$, the symbol $P_k(e)$ denotes the \emph{Peirce-$k$ projection} of $E$ onto $E_k(e)$. Peirce projections are contractive \cite[Corollary 1.2]{friedman85structure}. It is further known that Peirce subspaces are precisely the eigenspaces of the operator $L(e,e)$ and $L(e,e) = P_2(e) + \frac12 P_1(e).$ The Peirce-$2$ subspace $E_2(e)$ is a JB$^*$-algebra with product $x\circ_e y=\{x,e,y\}$ and involution $x^{*_e}=\{e,x,e\}$ (cf. \cite[Theorem 2.2]{braun78holomorphic} and \cite[Theorem 3.7]{kaup77jordan}). Triple products among elements in Peirce subspaces obey suitable rules known as \emph{Peirce arithmetic} described as follows: $$\{{E_{i}(e)},{E_{j}
(e)},{E_{k} (e)}\}\subseteq E_{i-j+k} (e),$$ if $i-j+k \in \{
0,1,2\}$ and is zero otherwise, and  $$\{{E_{2}
(e)},{E_{0}(e)},{E}\} = \{ {E_{0} (e)},{E_{2}(e)},{E}\} =0.$$


Any two tripotents $e, v\in E$ are called \emph{orthogonal} ($e\perp v$ in short) if $e\in E_0(v)$ or, equivalently, $v\in E_0(e)$ (equivalently, $L(e,v)=0$). We write $e\leq v$ if $v-e$ is a tripotent orthogonal to $e,$ equivalently, $e$ is a projection in the JB$^*$-algebra $E_2(v)$. If $e\perp v$, the sum $e+v$ is a tripotent satisfying $e,v\leq e+v$.


A non-zero tripotent $e$ in a JB$^*$-triple $E$ is called \emph{minimal} if $E_2(e)=\mathbb{C} e$, \emph{complete} or \emph{maximal} if $E_0(e)=\{0\}$, and \emph{unitary} if $E = E_2(e)$. By \cite[Proposition 3.5]{kaup77jordan}, the complete tripotents of a JB$^*$-triple $E$ coincide with the extreme points of its closed unit ball. That is, if we write $\partial_e\left(\mathcal{B}_{E}\right)$ for the extreme points of the closed unit ball of $E$, we have $\partial_e\left(\mathcal{B}_{E}\right) = \left\{\hbox{complete tripotents in } E \right\}.$

A JB$^*$-triple may not contain a single non-zero tripotent, but by the Krein-Milman theorem complete tripotents are ``abundant'' in every JBW$^*$-triple.


A closed subspace $I$ of a JB$^*$-triple $E$ is a \emph{triple ideal} if $\{E,E,I\}+\{E,I,E\} \subseteq I.$ It is known that $I$ is a triple ideal  if and only if $\{E,E,I\} \subseteq I$ (cf. \cite[Proposition 1.4]{dineen88centroid} and \cite[Proposition 1.3]{BuChu92}).  A closed subspace $I$ is called an \emph{inner ideal} of $ E$ if $\{I, E, I\} \subseteq I.$ For each tripotent $e$ in $E$, it can easily be deduced from Peirce arithmetic that $E_2(e)$ and $E_0(e)$ are inner ideals of $E$. The reader can consult \cite{EdRu92InnerIdealsTriples} for more details on inner ideals.

\section{Commutativity in JB$^*$-triples revisited}\label{sec: commutativity}

A JB$^*$-triple $E$ is called \emph{commutative} or \emph{abelian} or \emph{associative} if $$[L(a,b),L(x,y)]=L(a,b)L(x,y)-L(x,y)L(a,b)=0,$$ for all $a,b,x,y\in E,$ that is, $L(a,b)$ and $L(x,y)$ commute in the associative Banach algebra $B(E)$ (cf. \cite{kaup83riemann,hornIdeals, FriRuCommutative, dineen88centroid}).

The usual product in a JB$^*$-algebra $\A$ is always commutative in the sense that $a\circ b=b\circ a$ for every $a,b\in \A$. The JB$^*$-algebra $\A$ is called \emph{associative} if the any two Jordan multiplication operators $M_a, M_b:\A\rightarrow\A$ commute in $B(\A)$, i.e., $[M_a, M_b](x)=M_a M_b(x)-M_b M_a(x)=0$ for all $a,b,x\in\A$. It is well known that when a C$^*$-algebra $A$ is regarded as JB$^*$-algebra with respect to the natural Jordan product, $a\circ b =\frac12 (a b + ba)$, it is commutative as C$^*$-algebra if and only if it is associative as JB$^*$-algebra.  

The following remark will be employed on several occasions.

\begin{remark}{\rm \label{rmk: associative JB*-algebra commutative JB*-triple} Clearly each JB$^*$-subtriple of a commutative JB$^*$-triple is commutative. A JB$^*$-algebra $\A$ is associative if and only if it is commutative when regarded as a JB$^*$-triple {\rm(}see \cite[Proposition 4.1]{dineen88centroid}{\rm)}.}
\end{remark}

In the setting of C$^*$-algebras, a celebrated result due to Kaplansky proves that commutativity is characterized by the non-existence of non-zero $2$-nilpotent elements (cf. \cite[2.12.21]{diximer77C*-algebras}). In the case of JB$^*$-algebras (or in the wider setting of non-necessarily commutative JB$^*$-algebras), Kaplansky's characterization was obtained by Iochum, Loupias and Rodríguez-Palacios, who proved the following:

\begin{thm}[{\cite[Theorem 1]{iochum89commutativity}}] \label{thm: associative JB*-algebra nilpotent}
Let $\A$ be a JB$^*$-algebra. Then, $\A$ is associative if and only if for every element $a\in \A$ the condition $a^2=0$ implies $a=0$.
\end{thm}

The following lemma can be obtained by a straight application of the separate weak$^*$-continuity of the triple product in a JBW$^*$-triple.

\begin{lem}\label{l wstar dense subtriples} Let $E$ be a weak$^*$-dense JB$^*$-subtriple of a JBW$^*$-triple $M$. Then $E$ is commutative if and only if $M$ is.  
\end{lem}

\begin{remark}\label{r atomic decomp}We recall next another of the tools required in our arguments. Each JBW$^*$-triple $M$ decomposes as the direct sum ($M=\At\oplus^\infty\Nat$) of two (possibly zero) orthogonal weak$^*$-closed ideal $\At_{M} = \At$ and $\Nat_{M}=\Nat$ (called the \emph{atomic} and the \emph{non-atomic} part of $M$, respectively) such that $\At$ is precisely the weak$^*$-closed linear span of all minimal tripotents in $M$, $\Nat$ contains no minimal tripotents, the closed unit ball of $\Nat$ contains no extreme points, $\At_{*}$ coincides with norm closed linear span of all extreme points in the closed unit ball of $M_*$ and $\phi|_{\Nat} \equiv 0$ for all extreme point $\phi$ of the closed unit ball of $M_*$ (cf. \cite[Theorems 1 and 2]{friedman85structure}). 
Extreme points of the closed unit ball of $M_*$ are in one-to-one correspondence with minimal tripotents in $M$, actually for each $\phi\in \partial_e \left(\mathcal{B}_{M_*}\right)$ there exists a unique minimal tripotent $e =e_{\phi}$ such that \begin{equation}\label{eq Peirce minimal} P_2 (e) (x)  = \phi (x) e, \ \hbox{ for all } x\in M \hbox{ (see \cite[Proposition 4]{friedman85structure})}.
\end{equation} Furthermore, the atomic part decomposes as the orthogonal sum of a family $\{C_i\}_{i\in \Lambda}$ of Cartan factors, i.e., $\At=\bigoplus_{i\in\Lambda}^\infty C_i$ (see \cite[Proposition 2 and Theorem E]{friedman86gelfand}).

Given a JB$^*$-triple $E$, the above decomposition applies to its bidual, $E^{**},$  which is a JBW$^*$-triple. Let us note that in this case $\At_{E^{**}}\neq \{0\},$ because $\partial_e \left(\mathcal{B}_{E^*}\right) \neq \{0\}.$ Moreover, if $\pi_{\At}$ stands for the natural projection of $E^{**}$ onto its atomic part and $\iota_{E} : E\hookrightarrow E^{**}$ denotes the canonical embedding, the mapping $\Psi_{E}= \pi_{\At}\circ \iota_{E} : E\hookrightarrow \At_{E^{**}}$ is an isometric triple homomorphism with weak$^*$-dense image \cite[Proposition 1 and its proof]{friedman86gelfand}. The elements in $\partial_e \left(\mathcal{B}_{E^*}\right)$ are called \emph{pure atoms}. 
\end{remark}

By the Gelfand theory for JB$^*$-triples (see \cite[Corollary 1.11]{kaup83riemann}), each commutative JB$^*$-triple $E$ can be (isometrically) identified, via a triple isomorphism (i.e., a linear bijection preserving triple product), with the norm closed subspace of a $C_0(L)$ consisting of all $\mathbb{T}$-homogeneous (or $\mathbb{T}$-equivariant) continuous functions on a principal $\mathbb{T}$-bundle $L$, that is, a $\mathbb{T}$-symmetric (i.e., $\mathbb{T} L = L$) subset of a locally convex Hausdorff complex linear space $L$ such that $0 \notin L$ and $L \cup \{0\}$ is compact, i.e.,
\[E\cong C_0^\mathbb{T}(L):=\{a\in C_0(L):a(\lambda t)=\lambda a(t)\text{ for every } (\lambda,t)\in\mathbb{T}\times L\}.\]
The space $C_0^\mathbb{T}(L)$ is equipped with the supremum norm and the triple product given by $\{a,b,c\} = a \overline{b} c$,

Every commutative C$^*$-algebra is obviously a commutative JB$^*$-triple (i.e. a $C^{\mathbb{T}}_0 (L)$-space for an appropriate principal $\mathbb{T}$-bundle $L$, see \cite[Proposition 10]{Ol74} or \cite[Lemma 3.1]{FriedmanRusso82TAMS}). The example exhibited in \cite[Corollary 1.13 and subsequent comments]{kaup83riemann} shows that the class of commutative JB$^*$-triples is strictly wider.

The case of single generated JB$^*$-triples is more favourable since, for each element $x$ in a JB$^*$-triple $E$, the JB$^*$-subtriple generated by the element $x$, denoted by $E_x$, is (isometrically) JB$^*$-triple isomorphic to a commutative C$^*$-algebra $C_0
(\Omega)$ for some locally compact Hausdorff space $\Omega$
contained in $(0,\|x\|],$ such that $\Omega\cup \{0\}$ is compact,
where $C_0 (\Omega)$ denotes the Banach space of all
complex-valued continuous functions vanishing at $0.$ It is also
known that there exists a triple isomorphism $\Psi$ from $E_x$
onto $C_{0}(\Omega),$ satisfying $\Psi (x) (t) = t$ $(t\in \Omega)$ (cf.
\cite[Corollary 4.8]{Kaup77MathAnn}, \cite[Corollary 1.15]{kaup83riemann} and
\cite{FriedmanRusso82TAMS,FriRuCommutative}). The set $\overline{\Omega }=\hbox{Sp} (x)$ is
called the \emph{triple spectrum} of $x$. We should note that $C_0
(\hbox{Sp} (x)) = C(\hbox{Sp}(x))$, whenever $0\notin \hbox{Sp}
(x)$. As a consequence of this ``continuous triple functional calculus'', for each $x\in E$, there exists a unique element $y\in
E_x$ satisfying that $\{y,y,y\} =x$. The element $y$ will be denoted by
$x^{[\frac13 ]}$, and is termed the \emph{cubic root} of $x$. We can
inductively define $x^{[\frac{1}{3^n}]} =\left(x^{[\frac{1}{3^{n-1}}]}\right)^{[\frac{1}{3}]}$, $n\in \mathbb{N}$.
The sequence $(x^{[\frac{1}{3^n}]})$ converges in the weak$^*$
topology of $E^{**}$ to a tripotent denoted by $r(x)$ and called
the \emph{range tripotent} of $x$. The tripotent $r(x)$ is the
smallest tripotent $e\in E^{**}$ satisfying that $x$ is positive
in the JBW$^*$-algebra $E^{**}_{2} (e)$ (compare \cite[Lemma
3.3]{edwards96compact}).\smallskip

The following lemma gathers some facts from the folklore in JB$^*$-triples. It is stated here for completeness and to provide an explicit argument to the reader. 

\begin{lem}\label{l comm and unital} Let $E$ be a commutative JB$^*$-triple admitting a unitary tripotent $u$. Then $E$ is a commutative unital C$^*$-algebra. Consequently, each commutative JBW$^*$-triple is a commutative von Neumann algebra. 
\end{lem}

\begin{proof} It is easy to check from the commutativity of the JB$^*$-triple $E = E_2(u)$ that the Jordan product $\circ_{u}$ is associative, and hence the statement follows. For the second statement we observe that if $M$ is a commutative JBW$^*$-triple,   $\partial_e\left(\mathcal{B}_{M}\right)\neq \emptyset,$ and thus, we can find a complete tripotent $u$ in $M$. The representation $M$ as $C_0^\mathbb{T}(L)$ for some principal $\mathbb{T}$-bundle $L$, implies that $P_1(u) =0,$ and consequently $u$ is unital.
\end{proof}


We shall also handle in this note the inner ideal $E(a)$ generated by a single element $a$ in a JB$^*$-triple $E$. It is known that $E(a)$ is precisely the norm closure of the set $Q (a)(E) = \{a,E,a\}.$ Moreover $E(a)$ is a JB$^*$-subalgebra of $ E_2^{**}(r(a))$ and contains $a$ as a
positive element (cf. \cite[Proposition 2.1]{BuCHuZa2000MathScand}). In this case $E_a$ is nothing but the JB$^*$-subalgebra of $E(a)$ generated by the element $a$. Henceforth, we shall always regard $E(a)$ with its natural structure of JB$^*$-algebra described in this paragraph. Let us observe that if $e$ is a tripotent in $E$, the Peirce-2 subspace $E_2(e)$ is the inner ideal generated by $e$. Another interesting example is given by a complex Hilbert space $H$ regarded as a type 1 Cartan factor. It is known that the triple product in this case is given by $$\{\xi,\eta, \zeta\} =\frac12 \left( \langle \xi,\eta\rangle \zeta + \langle \zeta,\eta\rangle \xi \right), $$ where $\langle \cdot, \cdot\rangle $ denotes the inner product of $H$. Clearly, the set of non-zero tripotents in $H$ is the whole unit sphere with $H_2(\xi) = \mathbb{C}\xi$ and $H_1(\xi) = \{\xi\}_{\perp_2} = \{\eta\in H : \langle \eta,\xi\rangle =0\}$ for each norm-one element $\xi$ in $H$. It then follows that $H(\xi) = \mathbb{C} \xi$ for all $\xi\in H$. The inner ideal generated by a single element is too reduced in this case.  \smallskip

Orthogonality also makes sense as a relation between general elements in a JB$^*$-triple $E$. We say that $a,b\in E$ are orthogonal ($a\perp b$ in symbol) if $L(a,b)=0.$ The following reformulations of the fact $a\perp b$ can be consulted in \cite[Lemma 1]{BurFerGarMarPe2008}:
\begin{equation}
\label{eq ref orthogo}\begin{array}{ccc}
  \{a,a,b\} =0; & a \perp r(b); & r(a) \perp r(b); \\
  & & \\
  E^{**}_2(r(a)) \perp E^{**}_2(r(b));\ \ \ & r(a) \in E^{**}_0 (r(b));\ \ \  & a \in E^{**}_0 (r(b)); \\
  & & \\
  b \in E^{**}_0 (r(a)); & E_a \perp E_b & \{ b,b,a\}=0.
\end{array}
\end{equation}

Let $M$ be a JBW$^*$-triple. Since the triple product of $M$ is separately weak$^*$-continuous, $M_2(e)$ is a JBW$^*$-algebra. Given any $\varphi\in M_*$, there exists an unique tripotent $e=s(\varphi)\in M$ such that $\varphi=P_2(e)\circ \varphi$ and $\varphi|_{M_2(e)}$ is a faithful normal positive functional on $M_2(e)$ (see \cite[Proposition 2]{friedman85structure}). The tripotent $e=s(\varphi)$ is called the \emph{support tripotent} of $\varphi$.

In the next theorem we gather some known characterizations of commutativity in JB$^*$-triples borrowed from \cite{kaup83riemann,FriRuCommutative,dineen88centroid}, together with some new equivalent reformulations which will be later employed to estimate the numerical index of a JB$^*$-triple. 

\begin{thm} \label{thm: commutativity characterizations JB*-triples} Let $E$ be a JB$^*$-triple. Then, the following are equivalent:
\begin{enumerate}[$(i)$]
\item $E$ is commutative.
\item $E^{**}$ is commutative.
\item The atomic part of $E^{**}$, $\At_{E^{**}},$ is commutative, equivalently, $\At_{E^{**}}\cong \ell_\infty$.
\item $E\cong C_0^\mathbb{T}(L)$ for some principal $\mathbb{T}$-bundle $L$.
\item For each functional $\varphi\in E^{*}$ with support tripotent $e = s(\varphi)$ we have $P_1 (e)=0$.
\item For each functional $\phi\in \partial_e\left(\mathcal{B}_{E^{*}}\right)$ with support tripotent $e = s(\phi)$ we have $P_1 (e)=0$.
\item Each functional $\phi\in \partial_e\left(\mathcal{B}_{E^{*}}\right)$ is a triple homomorphism or a triple anti-homomorphism {\rm(}i.e. $\phi\{a,b,c\} = -\phi (a) \phi(b) \phi (c)${\rm)}.
\item For each minimal tripotent $e\in E^{**}$, $P_1(e)  =0$.

\item For each tripotent $v\in E^{**}$, $E_2^{**}(v)$ is an associative JBW$^*$-algebra {\rm(}or a commutative JBW$^*$-triple{\rm)} and $L(v,v) (E)\subseteq E_2^{**}(v)$.
\item For each complete tripotent $e\in E^{**}$, $E_2^{**}(e)$ is an associative JBW$^*$-algebra and $L(e,e) (E)\subseteq E_2^{**}(e)$.

\item For each complete tripotent $e\in E^{**}$, the JBW$^*$-algebra $E_2^{**}(e)$ does not contain non-zero $2$-nilpotent elements and $L(e,e) (E)\subseteq E_2^{**}(e)$.
    
\item For each $a$ in $E$, the inner ideal $E(a)$ is an associative JB$^*$-algebra (or a commutative JB$^*$-triple) and contains the image of the operator $L(a,a)$.

\item For each $a$ in $E$, the inner ideal $E(a)$ does not contain non-zero $2$-nilpotent elements as JB$^*$-algebra and $L(a,a) (E)\subseteq E(a)$.

\item For each $a$ in $E$ we have $L(a,a) (E)\subseteq E(a)$.
\end{enumerate}
\end{thm}

\begin{proof}  The equivalences $(i)\Leftrightarrow (ii)\Leftrightarrow (iii)$ follow from the fact that $E$ embeds isometrically and weak$^*$-densely as subtriple of $E^{**}$ and $\At_{E^{**}}$ (cf. Lemma~\ref{l wstar dense subtriples}).  The equivalence with $(iv)$ has been already commented before the theorem.

The equivalence of $(v),$ $(vi)$ and $(vii)$ with all previous statements is proved in \cite[Theorem 4.3]{dineen88centroid}, and the equivalence with $(viii)$ follows from \eqref{eq Peirce minimal}. Therefore, all statements from $(i)$ to $(viii)$ are equivalent. 

If $E^{**}$ is commutative, for each tripotent $v\in E^{**},$ the Peirce-2 space $E^{**}_2(v)$ is clearly an associative JBW$^*$-algebra since it is a JBW$^*$-subtriple of $E^{**}$. To show the second statement, we observe that if $E\cong C_0^{\mathbb{T}}(L)$ for some principal $\mathbb{T}$-bundle $L$, it easily follows that $E^{**}$ is a commutative von Neumann algebra, and hence $P_1 (v) =0,$ and hence $L(v,v) E^{**} = E^{**}_2(v) + E_1^{**} (v) = E_2^{**}(v)$ for all tripotent $v\in E^{**}.$ This shows that $(ix)$ is deduced from any the previous statements.   

The implication $(ix)\Rightarrow (x)$ is obvious.

$(x)\Rightarrow (ii)$ Since every element in $\partial_e \left(\mathcal{B}_{E^{**}}\right)$ is a complete tripotent \cite[Proposition 3.5]{kaup77jordan}, we can always find a complete tripotent $e\in E^{**}$. $E_2^{**}(e)$ is associative by assumptions. Moreover, since we also know that $L(e,e)=P_2(e)$ because $P_1(e) (E)\subseteq E_2^{**} (e).$ This shows that $E^{**} = E_2^{**} (e)$ is an associative JBW$^*$-algebra, equivalently, a commutative von Neumann algebra. 

The equivalences $(x)\Leftrightarrow (xi)$ and $(xii)\Leftrightarrow (xiii)$ follow from Theorem~\ref{thm: associative JB*-algebra nilpotent}. It is also easy to check that $(iv)$ implies $(xii).$ The implication $(xiii)\Rightarrow (xiv)$ is clear.

We shall finally prove that $(xiv)\Rightarrow (iii)$. Assume that $(xiv)$ holds. If the atomic part of $E^{**},$ $\At_{E^{**}}$, contains a Cartan factor $C$ with rank $\geq 2$, it follows, for example, from \cite[Lemma 3.10 and the discussion prior to it]{FerPeAdv2018} that $C$ (and hence $\At_{E^{**}}$ and $E^{**}$) contains a JB$^*$-subtriple isometrically isomorphic to $M_2 (\mathbb{C})$ or to $S_3(\mathbb{C}) = \left\{ \left(\begin{array}{cc}
\alpha & \beta \\
\beta & \delta \\
\end{array}\right) : \alpha, \beta, \delta\in \mathbb{C}\right\},$ moreover, we can find two orthogonal minimal tripotents $e_{11},e_{22}$ and a rank-2 tripotent $e_{12}$ in $E^{**}$ such that the JB$^*$-subtriple generated by them identifies isometrically with $M_2 (\mathbb{C})$ or with $S_3(\mathbb{C})$, $e_{11} \cong  \left(\begin{array}{cc}
1 & 0 \\
0 & 0 \\
\end{array}\right),$ $e_{22} \cong  \left(\begin{array}{cc}
0 & 0 \\
0 & 1 \\
\end{array}\right),$ and $e_{12} \cong  \left(\begin{array}{cc}
0 & 1 \\
1 & 0 \\
\end{array}\right)$.  By applying Kadison's transitivity theorem for JB$^*$-triples in \cite[Theorem 3.3]{bunce06saito-tomima-lusin} there exist norm-one elements $a,b\in E$ such that $a = e_{11} + P_0 (e_{11}+ e_{22}) (a) = e_{11} + a_0$ and $b = e_{12} + P_0(e_{12}) (b)= e_{12} + b_0.$ Let us comment how to get $a$. By \cite[Theorem 3.3$(c)$]{bunce06saito-tomima-lusin} applied to the orthogonal minimal tripotents $e_{11}$ and $e_{22},$ we can find orthogonal norm-one elements $a,c \in E$ such that $a = e_{11}+P_0(e_{11}) (a)$ and $c = e_{22}+P_0(e_{22}) (c)$. We deduce from $a\perp c$ that $P_0(e_{11}) (a) =P_0 (e_{11}+ e_{22}) (a) $ (cf. \eqref{eq ref orthogo}).

In this case, by orthogonality, $$\begin{aligned}
L(a,a) (b) &= L(e_{11},e_{11}) (b) + L(a_0, a_0) (b) \\
&=  L(e_{11},e_{11}) (e_{12}) + L(a_0,a_0) (b_0) = \frac12 e_{12} + L(a_0,a_0) (b_0) = \frac12 e_{12} + c_0,
\end{aligned}
$$ with $e_{12} \perp  c_0$. This is impossible since, by assumptions, $L(a,a) (E)\subseteq E(a)\subseteq E_2^{**} (r(a))$ and $r(a) = e_{11} + r(a_0)$ with $e_{11},e_{12}, e_{11}+e_{22} \perp r(a_0)$ and $P_2(r(a)) (e_{12}) = 0$ (note that $L(r(a),r(a)) (e_{12}) = L(e_{11},e_{11}) (e_{12}) + L(r(a_0),r(a_0)) (e_{12}) = \frac12 e_{12}$). We have therefore shown that $\At_{E^{**}}$ does not contain Cartan factors of rank $\geq 2$. 

Suppose next that $\At_{E^{**}}$ contains a non-trivial Cartan factor of rank-one $C$. It is well known that each JB$^*$-triple of rank-one is a complex Hilbert space equipped with its structure of type 1 Cartan factor (see, for example, \cite[Proposition 4.5]{BuChu92} and \cite{ChuIo90} or \cite[Theorem 2.3 and \S 3]{BeLoPeRod2004}). If dim$(C)\geq 2$, we can find an orthonormal set $\{\xi_1,\xi_2\}$ in the Hilbert space $C$. Applying the structure of type 1 Cartan factors, $\xi_1$ and $\xi_2$ are minimal and complete tripotents in $C$ with $\xi_i\in C_1(\xi_j) \subseteq E^{**}_1(\xi_j)$ for $i\neq j.$ The tripotents $\xi_1$ and $\xi_2$ are also minimal in $E^{**}$. By a new application of Kadison's transitivity theorem for JB$^*$-triples \cite[Theorem 3.3]{bunce06saito-tomima-lusin}, we can find two norm-one elements $a,b\in E$ such that $a = \xi_1 + P_0 (\xi_1) (a) = \xi_1 + a_0$ and $b = \xi_2 + P_0 (\xi_2) (b) = \xi_2 + b_0$. Observe that $ P_0 (\xi_1) (a)$ and $P_0 (\xi_2) (b)$ have zero component in $C$ because $\xi_1$ and $\xi_2$ are complete in $C$. The identities  
$L(r(a),r(a)) (\xi_2) = L(\xi_1,\xi_1) (\xi_2) + L(r(a_0),r(a_0)) (\xi_2) = \frac12 \xi_2,$ and  
$L(a,a) (b) = L(\xi_1,\xi_1) (\xi_2) + L(a_0,a_0) (b_0) = \frac12 \xi_2+ L(a_0,a_0) (b_0),$ with $L(a_0,a_0) (b_0)$ in the orthogonal complement of $C$ in $E^{**}$, prove that $L(a,a) (b)\notin E_2^{**} (r(a))\supseteq  E(a),$ which contradicts our assumptions.
 
Summarizing, all Cartan factors in $\At_{E^{**}}$ must coincide with $\mathbb{C},$ and hence $\At_{E^{**}}$ is commutative.
\end{proof}

The proof of the above theorem actually gives an argument to obtain the following corollary.

\begin{cor} \label{cor: nilpotent in non-commutative JB*-triple} Let $E$ be a JB$^*$-triple. The following assertions are equivalent:
\begin{enumerate}[$(a)$]
    \item $E$ is non-commutative.
    \item One of the following statements holds: \begin{enumerate}[$(1)$]
        \item There is an element $a$ in $E$ such that the inner ideal $E(a)$ contains a non-zero $2$-nilpotent element $b$ as JB$^*$-algebra {\rm(}i.e., $\{b,r(a),b\}=0${\rm)}, equivalently, $E(a)$ is non-associative.
        \item Every single generated inner ideal of $E$ is an associative JB$^*$-algebra and the atomic part of $E^{**}$ reduces to a $\ell_{\infty}$-sum of Hilbert spaces and at least one of them has dimension greater than or equal to $2$. 
    \end{enumerate} 
    \item One of the following statements holds: \begin{enumerate}[$(1)$]
        \item The atomic part of $E^{**}$ contains a Cartan factor with rank $\geq 2.$
        \item The atomic part of $E^{**}$ only contains rank-one Cartan factors and at least one of them is a Hilbert space of dimension greater than or equal to $2$. 
    \end{enumerate}
\end{enumerate}
\end{cor}

It is known that the unique rank-one Cartan factor contained in the atomic part of the bidual of a JB$^*$-algebra is $\mathbb{C}$ (cf. \cite[Proof of Corollary 3.2]{Dang1992} or  \cite[Proof of Corollary 3.4]{FerMarPe2004}). Then the next corollary trivially holds from the previous results.

\begin{cor} \label{cor: commutativity JB*-algebra} Let $\A$ be a JB$^*$-algebra. Then $\A$ is associative if and only if for each  $a$ in $E$ the inner ideal $E(a)$ is 
associative with respect to its natural structure of JB$^*$-algebra. 
\end{cor}

\begin{remark} In a general JB$^*$-triple $E$, the condition ``$E(a)$ is an associative JB$^*$-algebra for all $a\in E$'' does not imply that $E$ is commutative. For example, a Hilbert space $H,$ regarded as a type 1 Cartan factor, satisfies that $H(a) = \mathbb{C} a$ is a commutative JB$^*$-algebra, however, $H$ is not commutative. In this case we have a rank-one JB$^*$-triple, similar examples with bigger ranks can be obtained with $c_0$-sums of rank-one type 1 Cartan factors. 
\end{remark}

The statements $(i)$, $(ii)$, $(iii),$ $(v)$--$(xi)$ in Theorem~\ref{thm: commutativity characterizations JB*-triples} characterize commutativity of a JB$^*$-triple in terms of its dual and its predual. The new statements $(xii)$--$(xiv)$ refer to the intrinsic structure of $E.$ Assuming that $E$ is in fact a JBW$^*$-triple we have an abundant collection of tripotents in $E$ and we can somehow relax the hypothesis.

We introduce some notation first. Henceforth, we shall write $\mathcal{Q}_3$ for the JB$^*$-subalgebra of $M_2(\mathbb{C})$ of all matrices of the form $\begin{pmatrix}\alpha & \beta \\ \gamma & \alpha\end{pmatrix}$ with $\alpha,\beta,\gamma\in\mathbb{C}$. 

\begin{thm}\label{thm: commutativity characterizations JBW*-triples} Let $M$ be a JBW$^*$-triple. Then, the following are equivalent:
\begin{enumerate}[$(i)$]
\item $M$ is commutative.
\item $M$ is a commutative von Neumann algebra.
\item For each tripotent $e\in M$ we have $L(e,e) (M)\subseteq M_2 (e)$ {\rm(}equivalently, $P_1 (e) =0${\rm)}.
\end{enumerate}
\end{thm}

\begin{proof} The equivalence $(i)\Leftrightarrow (ii)$ follows from Lemma~\ref{l comm and unital}. 
$(ii)\Rightarrow (iii)$ is clear. 

$(iii)\Rightarrow (i)$ Let us take a complete tripotent $e$ in $M$. By the hypothesis, $L(e,e)(M) \subseteq M(e) = M_2(e)$, and hence $M = M_2(e)$ is a JBW$^*$-algebra. If $M=M_2(e)$ is non-associative as JBW$^*$-algebra, by \cite[Corollary 12]{iochum89commutativity}, $M_2(e)$ contains $\mathcal{Q}_3$ as JB$^*$-subalgebra. Consider the tripotents $v = \begin{pmatrix}0 & 0 \\ 1 & 0\end{pmatrix}$ and $w = \begin{pmatrix}1 & 0 \\ 0 & 1\end{pmatrix}$ in $\mathcal{Q}_3$. It is easy to check that $L(v,v) (w) = \frac12 w \in M_1 (v)$, and hence $L(v,v) (M) \nsubseteq M(v) = M_2(v)$, which contradicts our assumptions. 
\end{proof}

It should be commented that statement $(iii)$ in Theorem~\ref{thm: commutativity characterizations JBW*-triples} cannot be replaced by ``$P_1(e) =0$ for all complete tripotent $e\in M$''. For example $M = M_2(\mathbb{C})$ is a non-commutative JB$^*$-algebra which satisfies this property, since each complete tripotent in $M$ is unitary. Similar examples can be found in the setting of non-commutative finite JBW$^*$-algebras in the sense of \cite{HamKalPeFinite}.

\begin{remark}\label{r links between the two main hypotheses} In Theorem~\ref{thm: commutativity characterizations JB*-triples} and Corollary~\ref{cor: nilpotent in non-commutative JB*-triple}, there are two main properties linked up with the commutativity of a JB$^*$-triple $E,$ or with its rebuttal, namely, 
\begin{enumerate}[$(a)$]
\item There exists $a\in E$ such that $E(a)$ contains a non-zero $2$-nilpotent element;
\item There exists $b\in E$ such that $L(b,b)\nsubseteq E(b)$.
\end{enumerate}
Clearly, $(a)$ implies that $E$ is non-commutative (cf. Theorem~\ref{thm: associative JB*-algebra nilpotent}), and hence $(a)\Rightarrow (b)$ in light of Theorem~\ref{thm: commutativity characterizations JB*-triples}.

However the implication  $(b)\Rightarrow (a)$ does not always hold. A counterexample can be given by a Hilbert space of dimension $\geq 2$ regarded as a type 1 Cartan factor. 
\end{remark}

\section{The numerical index of a JB$^*$-triple}\label{sec:num index JBstartriples}

This section is devoted to compute the numerical index of a JB$^*$-triple, which seems to be an incognito problem until now. The numerical index is well determined in the case of C$^*$-algebras, JB$^*$-algebras, and non-necessarily commutative JB$^*$-algebras. More concretely, a C$^*$-algebra $A$ satisfies that $n(A) = 1$ if $A$ is commutative and $n(A)=\frac12$ otherwise \cite{HuruyaPAMS1977}, and the same conclusion holds when $A$ is replaced by a (unital) non-commutative Jordan $V$-algebra \cite[Theorem 26]{RodPal1980} or by a non-necessarily commutative JB$^*$-algebra \cite{iochum89commutativity,KadMorRodPal2001} (see also \cite[Proposition 3.5.44 and comments in \S 2.1.47 and page 422]{CabRodBookV1}). The results in this section show that every commutative JB$^*$-triple has numerical index one (see Lemma~\ref{l the numerical index of a commutative JB*-triple is 1}), while for a JBW$^*$-triple $M$ we show that $n(M)=1$ if $M$ is commutative and $\frac1e \leq n(M) \leq \frac12$ otherwise (cf. Theorem~\ref{t numerical index in JBW*-triples}).  In the general setting we establish that each (non-commutative) JB$^*$-triple $E$ whose second dual contains a Cartan factor of rank greater than or equal to $2$ in its atomic part satisfies $\frac1e \leq n(E) \leq \frac12$ (see Theorem \ref{t numerical index non-commutative JB*-triples}). 

We have already commented that for each Banach space $X$ and each $T\in B(X)$ we have $v_{B(X)}(T) = v_{B(X^*)}(T^*)$, and thus the inequality $n(X^*)\leq n(X)$ always holds \cite[Section 9 Corollary 6]{BonDunBookI}. The reader should be warned that the inequality just presented can be, in general, strict (cf. \cite{BoykoKadMar2007numIndexDual,MarJFA2008}).

Let us begin our study on the numerical index of a JB$^*$-triple with the case of a commutative JB$^*$-triple $E$. As we have already commented in the previous section, $E$ is not, in general, isometrically isomorphic to a commutative C$^*$-algebra. So, a direct application of the known references is not clear. However, this case can be easily derived from the available literature. 

\begin{lem}\label{l the numerical index of a commutative JB*-triple is 1} The numerical index of each commutative JB$^*$-triple $E$ is one.
\end{lem}

\begin{proof} It is known that $n(E^{**})\leq n(E^*) \leq n(E) \leq 1.$ Since $E^{**}$ is a commutative von Neumann algebra (cf. Lemma~\ref{l comm and unital}), and the numerical index of a C$^*$-algebra $A$ is $1$ or $\frac12$ depending if $A$ is commutative or not \cite{HuruyaPAMS1977}, we have $n(E^{**}) = n(E^*) = n(E) =1.$
\end{proof}

The following result is a well known consequence of Hahn--Banach theorem and will be applied (sometimes implicitly) several times along the rest of the document.

\begin{lem} \label{lem: numerical radius in subspace} Let $Y$ be a closed subspace of a Banach space $X$. If $T\in B(X)$ satisfies $T(Y)\subset Y$, then $V_{B(Y)}(T|_{Y})\subset V(T)$. In particular, $v_{B(Y)}(T|_{Y})\leq v(T)$.
\end{lem}

We continue with another technical tool, whose statement is probably known in the available sources like \cite{RodPal1980,iochum89commutativity,KadMorRodPal2001}. It is usually obtained via a celebrated result on algebra numerical range due to Bouldin \cite{bouldin71numericalII}. Here we isolate a precise computation of the numerical radius of the Jordan multiplication operator associated with a $2$-nilpotent element, the proof is included for completeness and to offer an alternative approach via Cauchy--Schwarz inequality, which somehow avoids the employment of Bouldin's result.  

\begin{lem} \label{lem: numerical radius of multiplication by nilpotent} Let $\A$ be JB$^*$-algebra and let $b$ be a $2$-nilpotent element in $\A$. Then $v(M_b)=\frac12\|b\|$.
\end{lem}

\begin{proof} We can clearly assume that $b$ is non-zero with $\|b\| =1.$ Under these hypotheses, $\A$ is non-associative (cf. Theorem~\ref{thm: associative JB*-algebra nilpotent}), and hence $n(\A)=\frac12$ by \cite[Theorem 5]{iochum89commutativity} (alternatively, \cite[Theorem 26]{RodPal1980}, \cite[Proposition 2.6]{KadMorRodPal2001} or \cite[Proposition 3.5.44 and comments in \S 2.1.47 and page 422]{CabRodBookV1}), thus $v_{_{B(\A)}}(M_b)\geq \frac12.$ Let us see the reciprocal inequality. 

Now, take any $\phi\in S_{\A^*}$ and $z\in S_\A$ such that $\phi(z)=1$. We have $\phi\circ M_b (z) = \phi({b}\circ z)=\varphi({b})$, where $\varphi$ is the functional on $\A$ given by $\varphi(y)=\phi(y\circ z)$ for all $y\in \A$. Since  $\|\varphi\|\leq\|\phi\|\|z\|\leq 1$ and $\varphi(\mathbf{1})=\phi(\mathbf{1}\circ z)=\phi(z)=1$, $\varphi$ is a state of $\A$ (see \cite[Lemma 1.2.2]{hanche-olsen84jordan}, if $\A$ is non-unital, we can take an approximate unit instead of $\mathbf{1}$ \cite[Lemma 1.2.2 and \S 3.6]{hanche-olsen84jordan}).
    
Let $\mathfrak{B}$ denote the JB$^*$-subalgebra of $\A$ generated by ${b}$. By the Shirshov--Cohn theorem \cite[Theorem 2.4.14]{hanche-olsen84jordan} (see also \cite[Corollary 2.2]{Wrig77}), $\mathfrak{B}$ is a JC$^*$-algebra, that is, a JB$^*$-subalgebra of some unital C$^*$-algebra $A$. Since $\varphi$ is a state of $\mathfrak{B}$, by \cite[Lemma 3.6.6]{hanche-olsen84jordan} we can find an extension of $\varphi$ to the unitization of $\mathfrak{B}$ inside $A$ whose norm is then  attained at the unit, the Hahn--Banach extension of the latter to $A$ nas norm-one and attains it at the unit element of $A$, therefore there is a state $\tilde{\varphi}$ of $A$ whose restriction to $\mathfrak{B}$ is $\varphi$.
    
We shall work now in the C$^*$-algebra $A$, where clearly ${b}^2=0$ and $\tilde{\varphi}$ is a state. Since ${b}$ is $2$-nilpotent, ${b} {b}^*\perp {b}^* {b}$ in $A$. Thus, the corresponding range projections $l=r({b} {b}^*)$ and $r=r({b}^* {b})$ in $A^{**}$ are orthogonal.  It is well known that ${b}=l {b} r$.

By the Cauchy--Schwarz inequality, we have
    \begin{align*}
        |\tilde{\varphi}({b})|^2&=|\tilde{\varphi}(l {b} r)|^2\leq\tilde{\varphi}(l)\tilde{\varphi}(r {b}^*  {b} r)\leq \tilde{\varphi}(l)\tilde{\varphi}(\|{b}^* {b}\|r)=\|{b}^* {b}\|\tilde{\varphi}(l)\tilde{\varphi}(r)=\|{b}\|^2\tilde{\varphi}(l)\tilde{\varphi}(r)\\&=\tilde{\varphi}(l)\tilde{\varphi}(r) \leq \tilde{\varphi}(l)\tilde{\varphi}(1-l)=\tilde{\varphi}(l)(1-\tilde{\varphi}(l))\leq 1/4,
    \end{align*}
    since $\tilde{\varphi}(l)\in[0,1]$.

    By combining all together, we have $|\phi\circ M_{b} (z)|=|\varphi({b})|=|\tilde{\varphi}({b})|\leq 1/2$. The arbitrariness of $\phi\in S_{\A^*}$ and $z\in S_\A$ with $\phi(z)=1$ assures that $v_{_{B(\A)}}(M_b)\leq \frac12$.    
\end{proof}

\begin{remark}\label{r numerical radious of multiplication operators} In the final part of the proof of the previous theorem we have rediscovered a conclusion close to Bouldin's results in \cite{bouldin71numericalII} via Cauchy--Schwarz inequality. Namely, if $b$ is a non-zero $2$-nilpotent element in a C$^*$-algebra $A$, the operators $L_b, R_b : A\to A$ given by $x\mapsto  b x, x b$ respectively both have numerical radius $\frac12 \|b\|.$  
\end{remark}

The next proposition is the key result to estimate  an upper bound of the numerical index of a non-commutative JBW$^*$-triple. 

\begin{prop}\label{prop: nilpotent numerical radius} Let $e$ be a nonzero tripotent in a JB$^*$-triple $E$. Suppose that there exists a non-zero $2$-nilpotent element $b\in E_2(e)$ {\rm(}i.e., $b^2\equiv\{b,e,b\}=0${\rm)}. Then, the operator $L(b,e):E\rightarrow E$ satisfies $v(L(b,e))=\frac 12 \|b\| = \frac12 \|L(b,e)\|$ in $B(E)$.
\end{prop}

\begin{proof} We may assume, without loss of generality, that $\|b\|=1$. It is shown in the Proof of \cite[Proposition 2.4]{bunce06saito-tomima-lusin} that there exist a unital JB$^*$-algebra $\A$ and an isometric triple embedding  $\Psi: E\hookrightarrow \A$ such that $\Psi(e)=p$ is a projection in $\A$. We shall denote $\tilde{b}=\Psi(b)\in \A_2(p)$. Clearly, $\tilde{b}$ is orthogonal to $1-p$. Observe that $v_{_{B(E)}}(L(b,e))=v_{_{B(\Psi(E))}}(L(\tilde{b},p))$.

By applying that $\Psi$ is a triple homomorphism we arrive to  $\tilde{b}^2=\{\tilde{b},1,\tilde{b}\}=\{\tilde{b},p+(1-p),\tilde{b}\}=\{\tilde{b},p,\tilde{b}\}=\Psi\{b,e,b\}=0$. In addition, the operator $L(\tilde{b},p)$ can be regarded as an operator in $B(\A)$, because $\Psi(E)$ is a JB$^*$-subtriple of $\A$. For all $a\in \A$, we get
    $$L(\tilde{b},p)(a)=\{\tilde{b},p,a\}=\{\tilde{b},1,a\}=\tilde{b}\circ a= M_{\tilde{b}} (a)$$
    because $1-p\perp b\in \A_2(p)$. That is, $L(\tilde{b},p)=L(\tilde{b},1) = M_{\tilde{b}}$ in $B(\A)$.

    Therefore, by applying Lemmata~\ref{lem: numerical radius in subspace} and  \ref{lem: numerical radius of multiplication by nilpotent} we deduce that $$v_{B(E)}\big(L(b,e)\big)=v_{B(\Psi(E))}\left(L(\tilde{b},p)\right)\leq v_{B(\A)}\left(L(\tilde{b},p)\right) = \frac12 \|M_{\tilde{b}}\|= \frac12 \|{\tilde{b}}\| = \frac12 .$$ 

Finally, since $b$ lies in the unital JB$^*$-algebra $E_2(e)$ ($e$ is the unit), we can write $b=h+ik$ with $h$ and $k$ self-adjoint in $E_2(e)$. Clearly, $\|b\|=1$ implies $\|h\|\geq 1/2$ or $\|k\|\geq 1/2$, we may suppose $\|h\|\geq 1/2$. We know that $\|h\|=\sup_{\phi\in \mathcal{S}(E_2(e))}|\phi(h)|$, where $\mathcal{S}(E_2(e))$ denotes the states of $E_2(e)$ (cf. \cite[Lemma 3.6.8]{hanche-olsen84jordan}). For any $\varepsilon>0$, we can find $\phi_\varepsilon\in \mathcal{S}(E_2(e))$ with $\phi_\varepsilon(e)=1$ such that $|\phi_\varepsilon(h)|>\frac12-\varepsilon$. Setting $\tilde{\phi}_\varepsilon=\phi_\varepsilon\circ P_2(e)\in S_{E^*}$, we have $\tilde{\phi}_\varepsilon(e)=1$ and $$\left|\tilde{\phi}_\varepsilon\circ L(b,e)(e)\right|=\left|\tilde{\phi}_\varepsilon\{b,e,e\}\right|=\left|\tilde{\phi}_\varepsilon(b)\right|=\left|\tilde{\phi}_\varepsilon(h)+i\tilde{\phi}_\varepsilon(k)\right|\geq \left|\tilde{\phi}_\varepsilon(h)\right|>\frac12-\varepsilon,$$
    which gives $v_{B(E)}(L(b,e))\geq 1/2$.
\end{proof}

The atomic decomposition of JBW$^*$-triples is an useful tool applied, for example, in the Gelfand--Naimark theorem for JB$^*$-triples \cite{friedman86gelfand}. There is, however, a finer decomposition of JBW$^*$-triples that mimics the lines of the famous Murray-von Neumann decomposition of von Neumann algebras, which is due to Horn and Neher \cite{horn87classification,HornNehClassCOnstJBWtriplesTAM1988}. Suppose we consider an algebraic property $(\mathcal{P})$ determining a concrete class of elements in a JBW$^*$-triple $M$. For example $(\mathcal{P})$: ``$e$ is a minimal tripotent'' or ``$e$ is an abelian tripotent, i.e. $M_2(e)$ is an associative JB$^*$-algebra''. Obviously, the elements in the class defined by the first property are inside the class given by the second one. By \cite[Corollary IV 3.61]{NeherBookGrid1987} (see also \cite[Theorem 4.5]{BattQJM1991}), there exist (unique) orthogonal w$^*$-closed ideals $I_{(\mathcal{P})}$ and $N_{(\mathcal{P})}$ such that $M = I_{(\mathcal{P})}\oplus ^{\infty} N_{(\mathcal{P})}$, $I_{(\mathcal{P})}$ coincides with the w$^*$-closure of the span of all elements of $M$ having
property $(\mathcal{P})$ and $N_{(\mathcal{P})}$ is the orthogonal complement of $I_{(\mathcal{P})}$. The class of all minimal tripotents gives rise to the atomic decomposition of $M,$ while the class of all abelian tripotents produces the decomposition of $M$ as the orthogonal sum of its \emph{type I part} and its \emph{continuous part}, that is, $M= M_{I}\oplus^{{\infty}} M_{c},$ where $M_{I}$ is the w$^*$-closure of the span of all abelian tripotents \cite{horn87classification} and, in case of being non-trivial, $M_{c}$ contains no non-zero abelian tripotents (cf. \cite[Proposition 4.13]{hornIdeals} or \cite[Corollary 4.11]{BattQJM1991}). $M$ is said to be of \emph{type I} (resp. \emph{continuous}) if $M=M_{I}$ (resp. $M= M_{c}\neq \{0\}$ to avoid ambiguity). This notation is perfectly compatible with the terminology in von Neumann algebra theory \cite{SakaiBook1971}, as well as in the classification of JBW$^*$-algebras \cite{hanche-olsen84jordan, AlfShultzBook2003}.

As the reader can already guess, the continuous part of a JBW$^*$-triple fits perfectly to the tools we have already developed in Proposition~\ref{prop: nilpotent numerical radius} without any further detail on its concrete form. However, in order to understand the type I part we shall need some more details. 

Let us consider two arbitrary von Neumann algebras $A\subset {B}(H)$ and $W\subset {B}(K)$, the algebraic tensor product $A\otimes W$ is canonically embedded into ${B}(H\otimes K)$, where $H\otimes K$ is the hilbertian tensor product of $H$ and $K$ (see \cite[Definition IV.1.2]{TakBookI}). The usual \emph{von Neumann tensor product} of $A$ and $W$ (denoted by $A\overline{\otimes} W,$) is precisely the von Neumann subalgebra of ${B}(H\otimes K)$ generated by the algebraic tensor product $A\otimes W$, that is, the weak$^*$ closure
of $A\otimes W$ in ${B}(H\otimes K)$ (see \cite[\S IV.5]{TakBookI}).\smallskip

The von Neumann tensor product also make sense for other types of operator spaces. Concretely, if $A$ denotes a commutative von Neumann algebra and $M$ is a JBW$^*$-subtriple of some $B(H)$, called a JW$^*$-triple (that is the case for Cartan factors of type 1, 2, 3, or 4, compare \cite{horn87classification} or \cite[\S 9]{HamKalPfitzPeJFA2020}). Following the standard notation \cite{horn87classification, HornNehClassCOnstJBWtriplesTAM1988}, we shall write
$A \overline{\otimes} M$ for the weak$^*$-closure of the algebraic tensor product $A\otimes M$ inside the usual von Neumann tensor product $A
\overline{\otimes} B(H)$ of $A$ and $B(H)$. Clearly $A \overline{\otimes} M$ is a JBW$^*$-subtriple of $A \overline{\otimes} B(H)$. 

By classification theory of type I JBW$^*$-triples established by Horn in \cite[Classification Theorem 1.7]{horn87classification}, the type I JBW$^*$-triples are precisely the
$\ell_{\infty}$-sums of JBW$^*$-triples of the form:
\begin{enumerate}[$(i):$]
\item $A \overline{\otimes} C$, where $A$ is an abelian von
Neumann algebra and $C$ is a Cartan factor realised as a
JW$^*$-subtriple of some $B(H)$;
\item $A \otimes C$ (algebraic tensor product) where $A$
is as before and $C$ is an exceptional Cartan factor.\smallskip
\end{enumerate}

Of course, $A \overline{\otimes} C = A {\otimes} C$ whenever $C$ is a finite dimensional non-exceptional Cartan
factor (see. \cite[Theorem IV.4.14]{TakBookI}). This particular setting, perhaps deserves an extra explanation. Let us take a finite-dimensional JB$^*$-triple $F$. To be coherent with the previos notation, we write $A \overline{\otimes} C$ for the injective tensor product of $A$ and $C$. If $A \cong C(\Omega)$ (with $\Omega$ the spectrum of $A$), the tensor product $A \overline{\otimes} C$ is identified with the space $C(\Omega, C)$, of all continuous functions on $\Omega$ with values in $C$ endowed with the pointwise operations and the supremum norm (cf. \cite[page 49]{RyanBook2013}). 

We shall frequently apply that given a unitary $u\in A$ and a tripotent $e$ in a Cartan factor $C$, the element $u\otimes e$ is a tripotent in $A \overline{\otimes} C$ and the corresponding Peirce subspaces satisfy \begin{equation}\label{eq Peirce decomp of elementary tensors}
\left( A \overline{\otimes} C\right)_j (u\otimes e) = A \overline{\otimes} \left(C_j (e)\right) 
\end{equation} (cf. \cite[Lemma 1.7$(i)$]{horn87classification}).

The Horn--Neher--type I--continuous classification has been shown a powerful tool to tackle problems in JB$^*$-triple theory, like the description of JB$^*$-triples satisfying the Dunford--Pettis and alternative Dunford--Pettis properties \cite{ChuMeDPPJLMS1997,AcPeQJM2001}, the study of the image of a contractive projection on a JBW$^*$-triple \cite{BuPeJMAA2002, NeRuTAMS2003, ChuNealRuJOT2004}, the study of measures of weak non-compactness in preduals of von Neumann algebras and JBW$^*$-triples \cite{HamKalPfitzPeJFA2020}, proof of Barton--Friedman conjecture for Grothendieck's inequalities in JB$^*$-triples \cite{HamKalPePfitzTAMS2021}, determination of finite JBW$^*$-algebras and triples \cite{HamKalPeFinite}, etcetera.  

It should be noted that, in most of cases, the von Neumann tensor product $A\overline{\otimes} C$ is out of scope for previous results computing the numerical index for $C(K,X)$ spaces \cite[Theorem 5]{MarPay2000Studia}, $L_{\infty}(\mu,X)$ spaces for a $\sigma$-finite measure $\mu$ \cite{MarVillLinfty2003}, and projective and injective tensor products of Banach spaces \cite{MarMerQuerIndexTensorProd2021}.

Before stating our estimation of the numerical index of a JBW$^*$-triple, we recall that a Banach space $X$ is called \emph{$L$-embedded} if its second dual writes as an $\ell_1$-sum of $X$ and some other closed subspace (see \cite{HarmWerWerBookMideals}). Obviously, reflexive Banach spaces are $L$-embedded. Preduals of von Neumann algebras are examples of $L$-embedded spaces \cite[Example IV.1.1]{HarmWerWerBookMideals}, and the same occurs to preduals of JBW$^*$-triples \cite[Proposition 3.4]{barton86weak} and preduals of real JBW$^*$-triples \cite[Proposition 2.2]{BeLoPeRod2004}. Martín proved in \cite[Theorem 2.1]{MartProcAMS2009Lembedded} that for each $L$-embedded space $X$, we have $n(X) = n(X^*)$. This holds, in particular, when $X$ is the predual of a (real or complex) JBW$^*$-triple.

\begin{thm}\label{t numerical index in JBW*-triples} Let $M$ be a JBW$^*$-triple. Then $n(M) = n(M_*) = 1$ if $M$ is commutative, and $\frac{1}{e}\leq n(M)= n(M_*)\leq \frac12$ otherwise.
\end{thm}

\begin{proof}  As we remarked in the comments preceding this theorem, the equality $n(M_*) = n(M)$ follows from \cite[Theorem 2.1]{MartProcAMS2009Lembedded}. If $M$ is commutative, we have $n(M)=1$ by Lemma~\ref{l the numerical index of a commutative JB*-triple is 1}. Let us write $M$ in the form \begin{equation}\label{eq typeI conts decomp of M} M = \left(\bigoplus^{\infty} A_j\overline{\otimes} C_j \right) \bigoplus^{\infty} M_{c},
\end{equation} where each $A_j$ is a commutative von Neumann algebra, each $C_j$ is a Cartan factor and $M_c$ is the continuous part of $M$ (cf. \cite{horn87classification,HornNehClassCOnstJBWtriplesTAM1988} and the paragraphs preceding this theorem). By \cite[Proposition 1]{MarPay2000Studia}, $n(M) = \inf\{n (A_j\overline{\otimes} C_j) : j\}\wedge n\left(M_c\right).$ If the continuous part $M_c$ is non-trivial we can clearly pick a non-zero tripotent $e\in M_c$ and the corresponding JB$^*$-algebra $\left(M_c\right)_2 (e)$ is non-associative, so by Theorem~\ref{thm: associative JB*-algebra nilpotent} there exists a non-zero $b\in \left(M_c\right)_2 (e)$ with $\{b,e,b\} =0.$ Proposition~\ref{prop: nilpotent numerical radius} asserts that the operator $L(b,e):M_c\rightarrow M_c$ satisfies $v(L(b,e))=\frac 12 \|b\| = \frac12 \|L(b,e)\|$ in $B(M_c)$ (and also in $B(M)$ when $L(b,e)$ is regarded as an operator in $B(M)$). We have therefore conclude that $\frac{1}{e}\leq n(M)\leq \frac12$ in this case. 

If in the decomposition \eqref{eq typeI conts decomp of M} there exists a Cartan factor $C_{j_0}$ with rank $\geq 2,$ $C_{j_0}$ contains a JB$^*$-subtriple isometrically isomorphic to $M_2 (\mathbb{C})$ or to $S_3(\mathbb{C}) = \left\{ \left(\begin{array}{cc}
\alpha & \beta \\
\beta & \delta \\
\end{array}\right) : \alpha, \beta, \delta\in \mathbb{C}\right\},$ and there exist two orthogonal minimal tripotents $e_{11},e_{22}$ and a rank-2 tripotent $e_{12}$ in $C_{j_0}$ such that the JB$^*$-subtriple of $C_{j_0}$ generated by them is JB$^*$-triple isometrically isomorphic to $M_2 (\mathbb{C})$ or to $S_3(\mathbb{C})$, $e_{11} \cong  \left(\begin{array}{cc}
1 & 0 \\
0 & 0 \\
\end{array}\right),$ $e_{22} \cong  \left(\begin{array}{cc}
0 & 0 \\
0 & 1 \\
\end{array}\right),$ and $e_{12} \cong  \left(\begin{array}{cc}
0 & 1 \\
1 & 0 \\
\end{array}\right)$ (cf. \cite[Lemma 3.10 and the preceding discussion]{FerPeAdv2018}). By fixing a unitary $u\in A_{j_0},$ the tripotents $\mathbf{e} = u\otimes e_{12}$ and $\mathbf{b} =  u\otimes e_{11}$ in $A_{j_0}\overline{\otimes} C_{j_0}$ satisfy $\mathbf{b} \in \left(A_{j_0}\overline{\otimes} C_{j_0}\right)_2 (\mathbf{e})$ with $\{\mathbf{b}, \mathbf{e}, \mathbf{b}\} =0$ (see \eqref{eq Peirce decomp of elementary tensors} or \cite[Lemma 1.7$(i)$]{horn87classification}). A new application of Proposition~\ref{prop: nilpotent numerical radius} implies that the operator $L(\mathbf{b},\mathbf{e}):A_{j_0}\overline{\otimes} C_{j_0}\rightarrow A_{j_0}\overline{\otimes} C_{j_0}$ satisfies $v(L(\mathbf{b},\mathbf{e}))=\frac 12 = \frac12 \|L(\mathbf{b},\mathbf{e})\|$ in $B(A_{j_0}\overline{\otimes} C_{j_0}),$ and hence $n(M) \leq n\left(A_{j_0}\overline{\otimes} C_{j_0}\right) \leq \frac12.$

We can finally assume that all Cartan factors in \eqref{eq typeI conts decomp of M} have rank-one. If one of them, lets say $C_{j_0}$ is a Hilbert space $H$ of dimension $\geq 2$ regarded as a type 1 Cartan factor. Consider a minimal orthogonal projection $p$ of $H$ onto a one-dimensional subspace, and let us write $\mathbf{1}$ for the unit of $A_{j_0}$. The element $\mathbf{p} :=\mathbf{1}\otimes p$ is a projection in the type I von Neumann algebra $A_{j_0}\overline{\otimes} B(H)$, and $A_{j_0}\overline{\otimes} C_{j_0} = A_{j_0}\overline{\otimes} H$ identifies isometrically with $\mathbf{p} \left(A_{j_0}\overline{\otimes} B(H)\right) \subset \left(A_{j_0}\overline{\otimes} B(H)\right)_2 (\mathbf{p}) \oplus \left(A_{j_0}\overline{\otimes} B(H)\right)_1 (\mathbf{p})$. Since dim$(H)\geq 2$, arguing as in the proof of Theorem~\ref{thm: commutativity characterizations JB*-triples}$(xiv)\Rightarrow (iii)$, we can find two minimal tripotents $\xi_1,\xi_2\in H$ with $\xi_1\in H_1(\xi_2)$ and $\xi_2\in H_1(\xi_1)$. We can further assume that $p=\xi_1\otimes\xi_1$. Consider the tripotents $\mathbf{e}_{11} = \mathbf{1}\otimes \xi_1$ and $\mathbf{e}_{12} = \mathbf{1}\otimes \xi_2$ in $A_{j_0}\overline{\otimes} H$ and in $A_{j_0}\overline{\otimes} B(H).$ We identify $\xi_j\in H=p B(H)$ with $\xi_1\otimes \xi_j$. We know that $\mathbf{e}_{11}$ is a projection in $A_{j_0}\overline{\otimes} B(H)$, $\mathbf{e}_{11}\in \left(A_{j_0}\overline{\otimes} B(H)\right) _1(\mathbf{e}_{12})$ and $\mathbf{e}_{12}\in \left(A_{j_0}\overline{\otimes} B(H)\right)_1(\mathbf{e}_{11})$. It is not hard to check that $\mathbf{e}_{12}^2 =0 $ in  $A_{j_0}\overline{\otimes} B(H)$ and $L(\mathbf{e}_{12}, \mathbf{e}_{11}) : A_{j_0}\overline{\otimes} B(H)\to A_{j_0}\overline{\otimes} B(H)$ is precisely the right multiplication operator by the $2$-nilpotent tripotent $\frac12 \mathbf{e}_{12}$, that is $L(\mathbf{e}_{12}, \mathbf{e}_{11}) (\mathbf{x}) = \frac12 \mathbf{x} \mathbf{e}_{12} = R_{\frac12 \mathbf{e}_{12}} ( \mathbf{x})$, for all $\mathbf{x}\in A_{j_0}\overline{\otimes} B(H).$ Having in mind Remark~\ref{r numerical radious of multiplication operators}, we conclude that $\frac14 = v_{_{B\left(A_{j_0}\overline{\otimes} B(H)\right)}} (R_{\frac12 \mathbf{e}_{12}}) = v_{_{B\left(A_{j_0}\overline{\otimes} B(H)\right)}} (L(\mathbf{e}_{12}, \mathbf{e}_{11}))$. Lemma~\ref{lem: numerical radius in subspace} implies that $v_{_{B\left(A_{j_0}\overline{\otimes} H\right)}} (L(\mathbf{e}_{12}, \mathbf{e}_{11})) \leq v_{_{B\left(A_{j_0}\overline{\otimes} B(H)\right)}} (L(\mathbf{e}_{12}, \mathbf{e}_{11})) = \frac14$. This proves that $n (M) \leq n\left( A_{j_0}\overline{\otimes} C_{j_0}\right)\leq \frac12.$

Finally, if $M_{c} =\{0\}$ and all Cartan factors in  \eqref{eq typeI conts decomp of M} coincide with $\mathbb{C},$ the JBW$^*$-triple $M$ is commutative and hence $n(M) =1.$
\end{proof}

\begin{remark}\label{r solution to question in MartMN2008} In \cite[page 384]{MarMathNach2008}, Martín asked whether in case that $Y$ is a JB$^*$-triple or the predual of a JBW$^*$-triple, the condition $n(Y) =1$ implies $n(Y^*) =1$. The result by Martín in \cite[Theorem 2.1]{MartProcAMS2009Lembedded} gives a positive solution in case that $Y$ is a JBW$^*$-triple predual. Our Theorem~\ref{t numerical index in JBW*-triples} adds a positive answer when $Y$ is a JBW$^*$-triple. 
\end{remark}

In the next corollary we extend \cite[Proposition 3.3]{MarMathNach2008} to the setting of JBW$^*$-triples.

\begin{cor}\label{c index 1 in preduals} Let $E$ be a JB$^*$-triple. Then, the following are equivalent:
\begin{enumerate}[$(i)$]
\item $n(E^*) = 1$.
\item $E^*$ has the alternative Daugavet property.
\item $E^{**}$ has the alternative Daugavet property.
\item $|\phi (u)| = 1$ for every $u \in \partial_e (\mathcal{B}_{E^{**}})$ and every $ \phi \in \partial_e (\mathcal{B}_{E^{*}})$.
\item The atomic part of $E^{**}$ is isometrically isomorphic to a commutative von Neumann algebra (obviously atomic).
\item $E$ (equivalently, $E^{**}$) is commutative.
\end{enumerate}

\end{cor}

\begin{proof} Since $E^{**}$ is a JBW$^*$-triple with predual $E^*,$ Theorem 2.3 in \cite{MarMathNach2008} assures the equivalence between statements $(ii), (iii), (iv)$ and $(v)$. The equivalence $(i)\Leftrightarrow (vi)$ follows from Theorem~\ref{t numerical index in JBW*-triples}. The implication $(i)\Rightarrow (ii)$ is a consequence of \cite[Lemma 2.3]{MarOikh2004}. (The implication $(vi)\Rightarrow (iv)$ can be also deduced from Theorem~\ref{thm: commutativity characterizations JB*-triples} and standard theory of commutative von Neumann algebras).     
\end{proof}

Proposition~\ref{prop: nilpotent numerical radius} has been a powerful tool to characterize commutativity of JBW$^*$-triples in terms of numerical index. However, the existence of tripotents in a general JB$^*$-triple is simple hopeless. In order to throw some new light to the general setting, our next goal will be the following strengthened version of Proposition~\ref{prop: nilpotent numerical radius}.

\begin{prop}\label{prop for a JBtriple admitting a non-commutative inner ideal} Let $E$ be a JB$^*$-triple. Suppose there exists an element $a$ in $E$ such that the inner ideal $E(a)$ contains a norm-one $2$-nilpotent element $b$ with respect to its natural structure of JB$^*$-algebra, that is, $\{b,r(a),b\}=0$. Let us consider the element $c := b\circ_{r(a)} b^{*_{r(a)}}$ which is non-zero and positive in the JB$^*$-algebra $E(a).$ Then the operator $T = L(b,c) + L(c,b^{*_{r(a)}})$ satisfies $\|T\|_{B(E)} = 1$ and $v_{_{B(E)}} (T) = \frac12.$ Consequently, $\frac1e \leq n(E)\leq \frac12$ in this case.  
\end{prop}

Before dealing with the proof of this proposition we shall consider a technical lemma. 

\begin{lem} \label{lem: b cube $2$-nilpotent} Let $b$ be a non-zero $2$-nilpotent  element in a JB$^*$-algebra $\A$. Then, the element $(b\circ b^*) \circ b=\frac{1}{2}\{b,b,b\}$ is non-zero and $2$-nilpotent. 
\end{lem}

\begin{proof} Let us begin from the first equality. Since the triple product of $\A$ is uniquely given by the expression in \eqref{eq JB*-algebra triple product} (cf. \cite[Proposition 5.5]{kaup83riemann}), it is easy to check from the assumptions on $b$ that $$\{b,b,b\}=2(b\circ b^*)\circ b-(b\circ b)\circ b^*=2(b\circ b^*)\circ b.$$ Clearly, $\{b,b,b\}$ is non-zero by the third axiom in the definition of JB$^*$-triples. 
 
By the Shirshov--Cohn theorem (see \cite[Theorem 2.4.14]{hanche-olsen84jordan} or  \cite[Corollary 2.2]{Wrig77}), the JB$^*$-subalgebra $\mathfrak{B}$ of $\A$ generated by $\{b,b^*\}$ is a JC$^*$-algebra, that is a JB$^*$-subalgebra of a C$^*$-algebra. Therefore, there exists a $C^*$-algebra $A$ such that $\mathfrak{B}$ is a JC$^*$-subalgebra of $A$. In this case, $b^2=0$ also in $A$, and
    \begin{align*}
        (b\circ b^*)\circ b=\left(\frac{bb^*+b^*b}{2}\right)\circ b=\frac{1}{4}(bb^*b+b^*bb+bbb^*+bb^*b)=\frac{1}{2}bb^*b=\frac{1}{2}\{b,b,b\}\in \A,
    \end{align*}
which implies that  
    $\big[(b\circ b^*)\circ b\big]^2=\frac{1}{4}(\{b,b,b\})^2=\frac{1}{4}(bb^*b)^2=\frac{1}{4}bb^*bbb^*b=0$ in $A,$ $\mathfrak{B}$ and $\A$.
\end{proof}

\begin{proof}[Proof of Proposition~\ref{prop for a JBtriple admitting a non-commutative inner ideal}]

We can clearly assume that $a$ is a norm-one element. As in the statement of the proposition, we set $c=b\circ_{r(a)} b^{*_{r(a)}}=\{b,b,r(a)\}\in E(a)\subset E$. To simplify the notation we shall write $\circ$ and $^*$ for the Jordan product $\circ_{r(a)}$ and the involution $*_{r(a)}$ of the JB$^*$-algebra $E(a),$ respectively.

Let us observe that, by Kaup's Banach--Stone theorem \cite[Proposition 5.5]{kaup83riemann}, the triple product among elements in the JB$^*$-algebra is uniquely determined by the expression in \eqref{eq JB*-algebra triple product}. By Lemma~\ref{lem: b cube $2$-nilpotent}, $\{b,b,b\}=b^{[3]}$ is a $2$-nilpotent element in $E(a)$.

We claim that \begin{equation}\label{eq b cube is a nilpotent }\hbox{$2 b\circ c =\{b,b,b\}=b^{[3]}$ also is a $2$-nilpotent element in the JB$^*$-algebra $E(b\circ b^*) = E(c)$}
\end{equation} (the reader should be warned that we have changed the Jordan product and the involution). Namely, by \cite[Theorem 10]{iochum89commutativity}, the JB$^*$-subalgebra $\mathfrak{B}$ of $E(a)$ generated by $b$ is JB$^*$-algebra (isometrically) isomorphic to $C_0(K,\mathcal{Q}_3)$ for some compact set $K\subset [0,\|b\|]$ containing $0$, where $\mathcal{Q}_3=\left\{\begin{pmatrix}\alpha & \beta \\ \gamma & \alpha\end{pmatrix}  : \alpha,\beta,\gamma\in\mathbb{C}\right\}.$ Furthermore, under this representation,  the element $b$ corresponds to the function $t\mapsto\begin{pmatrix}0 & 0 \\ t & 0\end{pmatrix}$. Thus, $b\circ b^*\equiv \begin{pmatrix}0 & 0 \\ t & 0\end{pmatrix}\circ \begin{pmatrix}0 & t \\ 0 & 0\end{pmatrix}= \begin{pmatrix}t^2/2 & 0 \\ 0 & t^2/2\end{pmatrix}$
and \begin{equation}\label{eq b circ c} 
    b^{[3]}=2(b\circ b^*)\circ b =2 c \circ b \equiv\begin{pmatrix}t^2 & 0 \\ 0 & t^2\end{pmatrix}\circ \begin{pmatrix}0 & 0 \\ t & 0\end{pmatrix}=\begin{pmatrix}0 & 0 \\ t^3 & 0\end{pmatrix}
\end{equation} Moreover, $r(c) = r(b\circ b^*)\equiv\begin{pmatrix}1 & 0 \\ 0 & 1\end{pmatrix}$ in $\mathfrak{B}^{**}$.

Since $\mathfrak{B}(b\circ b^*)=\mathfrak{B}^{**}(r(b\circ b^*))\cap\mathfrak{B}$, it contains $b^{[3]}$. Thus, $b^{[3]}$ lies in $E(b\circ b^*)$. We can compute in $\mathfrak{B}$ the operation
$$b^{[3]}\circ_{r(b\circ b^*)} b^{[3]}\equiv\left\{\begin{pmatrix}0 & 0 \\ t^3 & 0\end{pmatrix},\begin{pmatrix}1 & 0 \\ 0 & 1\end{pmatrix},\begin{pmatrix}0 & 0 \\ t^3 & 0\end{pmatrix}\right\}=\begin{pmatrix}0 & 0 \\ t^3 & 0\end{pmatrix}^2=\begin{pmatrix}0 & 0 \\ 0 & 0\end{pmatrix}$$
to conclude that $b^{[3]}$ is $2$-nilpotent in $\mathfrak{B}(b\circ b^*)$ and also in $E(b\circ b^*) =E(c)$. This concludes the proof of the claim in \eqref{eq b cube is a nilpotent }.

On the other hand, the element $c=b\circ_{r(a)} b^{*_{r(a)}}$ is positive in the JB$^*$-algebra $E(a)$ by definition, and since $E(a)$ is an inner ideal of $E$, we conclude that $Q(c)(E)=\{c,E,c\}\subset E(a)$, so $E(c)\subset E(a)$, and $r(c)$ coincides with the support projection in $E(a)^{**}$. Hence, $r(c)\leq r(a)$ in $E(a)^{**}$ (and in $E^{**}$), because $c$ is positive in the latter JB$^*$-algebra. Consequently, $E(c)$ is a JB$^*$-subalgebra of $E(a)$.

We obviously have $b\in \mathfrak{B}(b\circ b^*)\subset E(c)$. Thus, $b\circ_{r(c)}c=b\circ_{r(a)}c=b\circ_{r(a)}(b\circ_{r(a)}b^{*_{r(a)}})=\frac{1}{2}b^{[3]}$ is nilpotent in $\mathfrak{B}(b\circ b^*)$ ($\subset E(c)\subset E(a)$) by what we proved two paragraphs above. If follows that $b\circ_{r(c)}c$ is nilpotent in $E(c)\subset E(a)$.

As it is essentially established in the proofs of \cite[Proposition 2.4]{bunce06saito-tomima-lusin} and \cite[Lemma 2.5]{FerPeJLMS2006}, and explicitly in \cite[Lemma 3.9$(i)$]{EdFerHosPe2010}, there exists an isometric triple embedding $\Psi$ from $E^{**}$ into a unital JBW$^*$-algebra $\A$ mapping $r(c)$ to a projection $p$ in $\A.$ The Jordan product and involution of $\A$ will be denoted by $\circ$ and $*$, respectively. Observe that, by composing $\Psi$ with the natural embedding of $E$ into its bidual, we can also see $\Psi$ as an isometric triple embedding from $E$ into $\A.$ Since $\Psi (r(c))$ is a projection in $\A$ and $c$ is positive in $E^{**}_2 (r(c))$, it follows that $\Psi(c)=\tilde{c}\geq 0$ in $\A$. 

We turn now our attention to the bounded linear operator $T$ on $E^{**}$ defined by $T= L(b,c)+L(c,b^{*_{r(c)}})$. Since $b,c,$ and $b^{*_{r(c)}}$ all lie in $E$, the mapping $T$ maps $E$ into $E$ and can be regarded as an operator in $B(E).$ Observe that, by the separate weak$^*$-continuity of the triple product and Goldstine theorem, $(T|_E)^{**}=T$. Clearly, $\Psi|_{E} : E\to \Psi(E)$ is a surjective linear isometry and a triple isomorphism, and we can identify $B(E)$ and $B(\Psi(E))$ via the assignment $S\mapsto \tilde{S} :=\Psi S \Psi^{-1}$. Clearly, $\tilde{T} =\Psi T \Psi^{-1} = L(\tilde{b},\tilde{c})+L(\tilde{c},\tilde{b}^{*})$ with $\tilde{b}=\Psi(b)$, $\tilde{c} = \Psi(c)$, and the fact that ${*_{r(\tilde{c})}}$ coincides with the involution in $\A$ because $\tilde{c}$ is positive in $\A$ (recall that $\Psi (r(c)) =p$ is a projection in $\A$).  

Furthermore, since $\Psi(E)$ is a JB$^*$-subtriple of $\A$ and  $\tilde{b},\tilde{c},\tilde{c}^*\in \Psi (E),$ the mapping $\tilde{T} = L(\tilde{b},\tilde{c})+L(\tilde{c},\tilde{b}^{*})$ can be also regarded as an operator in $B(\A)$. According to the latter and having in mind the expression of the triple product in a JB$^*$-algebra in \eqref{eq JB*-algebra triple product}, $\tilde{T}\in B(\A)$ is of the form
\begin{align*}
    \tilde{T}(z)&=L(\tilde{b},\tilde{c})(z)+L(\tilde{c},\tilde{b}^{*_{r(\tilde{c})}})(z)=\{\tilde{b},\tilde{c},z\}+\{\tilde{c},\tilde{b}^{*_{r(\tilde{c})}},z\} \\    &=(\tilde{b}\circ\tilde{c})\circ z+(z\circ\tilde{c})\circ \tilde{b}-(\tilde{b}\circ z)\circ \tilde{c}+(\tilde{c}\circ\tilde{b})\circ z+(z\circ\tilde{b})\circ \tilde{c}-(\tilde{c}\circ z)\circ \tilde{b}\\
    &= 2 (\tilde{b}\circ\tilde{c})\circ z=2 M_{ \tilde{b}\circ\tilde{c}}(z) =M_{2 \tilde{b}\circ\tilde{c}}(z),
\end{align*}
for all $z\in \A$. Clearly, $\|\tilde{T}\|_{B(\A)}=\|M_{2 \tilde{b}\circ\tilde{c}}\|_{B(\A)}=\|2 \tilde{b}\circ\tilde{c}\|=\| \tilde{b}^{[3]}\|= \|\tilde{b}\|^3 = \|b\|^3=1$. Moreover, $$\|\tilde{T}|_{\Psi(E)}\|_{B(\Psi(E))} = \|T\|_{B(E)} \geq \|T(r(c))\|= \|\{b,c,r(c)\}+\{c,{b}^{*_{r({c})}},r(c)\}\|=\| 2 b\circ_{r(c)}c\| =1.$$
Since $\|T\|_{B(E)} = \|\tilde{T}|_{\Psi(E)}\|_{B(\Psi(E))} \leq \|\tilde{T}\|_{B(\A)} =1$, we deduce that 
$$\|T\|_{B(E)}=\|\tilde{T}|_{\Psi(E)}\|_{B(\Psi(E))}=\|\tilde{T}\|_{B(\A)}=1.$$

Now, by applying that $b\circ_{r(c)}c = \{b,r(c),c\}$ is a $2$-nilpotent element in $E(c)$ (see \eqref{eq b cube is a nilpotent }), we obtain that $\tilde{c}\circ_p \tilde{b} = \{\tilde{c}, p, \tilde{b} \} = \{\Psi({c}), \Psi(r(c)), \Psi({b}) \}$ is $2$-nilpotent in $\Psi(E)(\tilde{c})$, and hence in $\A(\tilde{c})$ and in $\A$, because $\Psi(r(c)) =p\geq r(\tilde{c})$ is a projection in $\A^{**}$. Next, we call to Lemmata~\ref{lem: numerical radius in subspace} and \ref{lem: numerical radius of multiplication by nilpotent} to obtain $v_{B(E)}(T) = v_{B(\Psi(E))}(\tilde{T}|_{\Psi(E)}) \leq v_{B(\A)}(\tilde{T})=v_{B(\A)}(M_{2 \tilde{b}\circ\tilde{c}})= \frac12\|2 \tilde{b}\circ\tilde{c}\|=\frac12=\frac12\|T\|_{B(E)}$. This shows that $n(E)\leq \frac12$.

Actually, the operator $T$ has numerical radius exactly $\frac12=\frac12\|T\|_{B(E)}$. Namely, since $2 b\circ c$ is $2$-nilpotent in $E(c)$ (cf. \eqref{eq b cube is a nilpotent }) and $T|_{E(c)} : E(c)\to E(c)$ coincides with the Jordan multiplication operator by the element $2 b\circ c,$ that is, $T(z) = 2(b\circ c)\circ z$, for all $z\in E(c),$ Lemmata~\ref{lem: numerical radius in subspace} and \ref{lem: numerical radius of multiplication by nilpotent} imply that $v_{B(E)}(T)\geq v_{B(E(c))}(T|_{E(c)}) =v_{B(E(c))}(M_{2 b\circ c})= \frac12\|M_{2 b\circ c}\|= \frac12=\frac12\|T|_{E(c)}\|=\frac12\|T\|$.
\end{proof}

We recall that an element $a$ in a JB$^*$-triple $E$ is called \emph{abelian} or \emph{commutative} if the inner ideal of $E$ generated by $a$, $E(a)$, is an associative JB$^*$-algebra (cf. \cite{horn87classification} and  \cite[page 196]{BunChuZalMathZ2002}). We say that $E$ is \emph{antiliminal} if it contains no non-zero abelian elements. This notion was introduced by Bunce, Chu and Zalar in their study of those JB$^*$-triples $E$ whose second dual is a type I JBW$^*$-triple in a clear analogy with the notions of liminal and postliminal C$^*$- and JB$^*$-algebras \cite[\S 3]{BunChuZalMathZ2002}.

\begin{cor}\label{c antiliminal JB*-triples have index 1/2} Let $E$ be a  JB$^*$-triple containing a non-commutative element. Then $\frac1e\leq n(E)\leq \frac12.$ The conclusion clearly holds if $E$ is antiliminal.
\end{cor}

\begin{proof} Suppose $E$ contains a non-commutative element $a.$ The inner ideal $E(a)$ is a non-asociative JB$^*$-algebra, and hence it contains a non-zero 2-nilpotent element (cf. Theorem~\ref{thm: associative JB*-algebra nilpotent}). The desired conclusion follows from Proposition~\ref{prop for a JBtriple admitting a non-commutative inner ideal}.     
\end{proof}

As we already commented in Theorem~\ref{thm: commutativity characterizations JB*-triples}, the commutativity of a JB$^*$-triple $E$ is violated by the existence of a non-commutative element in $E$. However, the pathology ``being non-commutative'' is actually characterized by the existence of an element $a\in E$ satisfying that $L(a,a)\nsubseteq E(a)$. We have computed the numerical index of a non-commutative JB$^*$-triple admitting a non-commutative element. We shall see next that we can improve a bit the conclusion.

\begin{thm}\label{t numerical index non-commutative JB*-triples} Let $E$ be a JB$^*$-triple then the following statements hold:
\begin{enumerate}[$(a)$]
\item $n(E) =1$ if $E$ is commutative. 
\item $n(E^*)=n(E^{**}) \leq \frac12$ if and only if $E$ is non-commutative.
\item If $E$ is non-commutative due to the presence of a Cartan factor of rank greater than or equal to $2$ in the atomic part of $E^{**}$, we have $\frac1e \leq n(E)\leq \frac12$.\end{enumerate}
\end{thm}

\begin{proof} After Lemma~\ref{l the numerical index of a commutative JB*-triple is 1}, we can always suppose that $E$ is non-commutative. Since $E$ is commutative if and only if $E^{**}$ is, the conclusion in $(b)$ follows from Theorem~\ref{t numerical index in JBW*-triples}.

In view of the results in section~\ref{sec: commutativity}, we can assume that one of the next statements holds:
\begin{enumerate}[$(1)$]
\item The atomic part, $\At_{E^{**}}$, of $E^{**}$ contains a Cartan factor of rank greater than or equal to $2$.
\item The atomic part of $E^{**}$ is an $\ell_{\infty}$-sum of complex Hilbert spaces regarded as type 1 Cartan factors and at least one of them has dimension greater than or equal to $2$. 
\end{enumerate}

We claim that assuming $(1)$, the JB$^*$-triple $E$ contains a non-abelian element, and hence the conclusion in $(c)$ will follow from Proposition~\ref{prop for a JBtriple admitting a non-commutative inner ideal}. So, we assume that $E^{**}$ contains a Cartan factor $C$ of rank $\geq 2$ in its atomic part. Arguing as in the proof of Theorem~\ref{thm: commutativity characterizations JB*-triples}~$(xiv)\Rightarrow(iii)$, we deduce that $C$ contains a JB$^*$-subtriple isometrically isomorphic to $M_2 (\mathbb{C})$ or to $S_3(\mathbb{C}) = \left\{ \left(\begin{array}{cc}
\alpha & \beta \\
\beta & \delta \\
\end{array}\right) : \alpha, \beta, \delta\in \mathbb{C}\right\},$ moreover, we can find two orthogonal minimal tripotents $e_{11},e_{22}$ and a rank-2 tripotent $e_{12}$ in $E^{**}$ such that the JB$^*$-subtriple generated by them identifies isometrically with $M_2 (\mathbb{C})$ or with $S_3(\mathbb{C})$ under a triple isomorphism satisfying $e_{11} \cong  \left(\begin{array}{cc}
1 & 0 \\
0 & 0 \\
\end{array}\right),$ $e_{22} \cong  \left(\begin{array}{cc}
0 & 0 \\
0 & 1 \\
\end{array}\right),$ and $e_{12} \cong  \left(\begin{array}{cc}
0 & 1 \\
1 & 0 \\
\end{array}\right)$ (cf. \cite[Lemma 3.10 and the discussion prior to it]{FerPeAdv2018}). We also justified in the quoted proof how, by Kadison's transitivity theorem for JB$^*$-triples in \cite[Theorem 3.3]{bunce06saito-tomima-lusin}, we can find norm-one elements $a,b\in E$ such that $a = e_{11} + P_0 (e_{11}+ e_{22}) (a) = e_{11} + a_0$ and $b = e_{12} + P_0(e_{12}) (b)= e_{12} + b_0.$ We observe that $E^{**}_{0} (e_{12}) = E_0^{**} (e_{11}+e_{22})$. The element $c = Q(b) (a) = e_{22} + Q(b_0) (a_0) = e_{22} + c_0$ lies in $E(b)$ with $c_0\in E^{**}_0(e_{12})$. Having in mind that $r(b) = e_{12} + r(b_0)$, it is not hard to check, by orthogonality, that $$\begin{aligned}
(c\circ_{r(b)} c^{*_{r(b)}}) \circ_{r(b)} c & = (e_{22}\circ_{e_{12}} e_{11}) \circ_{e_{12}} e_{22} +  (c_0\circ_{r(b_0)} c_0^{*_{r(b_0)}}) \circ_{r(b_0)} c_0 \\
&= \frac12 e_{12} \circ_{e_{12}} e_{22} +  (c_0\circ_{r(b_0)} c_0^{*_{r(b_0)}}) \circ_{r(b_0)} c_0 = \frac12 e_{22} +  (c_0\circ_{r(b_0)} c_0^{*_{r(b_0)}}) \circ_{r(b_0)} c_0,
\end{aligned}$$ 
and 
$$\begin{aligned}
(c\circ_{r(b)} c) \circ_{r(b)} c^{*_{r(b)}} & = (e_{22}\circ_{e_{12}} e_{22}) \circ_{e_{12}} e_{11} +  (c_0\circ_{r(b_0)} c_0) \circ_{r(b_0)} c_0^{*_{r(b_0)}} \\
&= 0 +  (c_0\circ_{r(b_0)} c_0) \circ_{r(b_0)} c_0^{*_{r(b_0)}} .
\end{aligned}$$ We deuce from the last two identities that $(c\circ_{r(b)} c^{*_{r(b)}}) \circ_{r(b)} c \neq (c\circ_{r(b)} c) \circ_{r(b)} c^{*_{r(b)}},$ and hence $E(b)$ is non-associative as JB$^*$-algebra. This concludes the proof of the theorem. 
\end{proof}

With our current knowledge and tools, we cannot conclude that every non-commutative JB$^*$-triple $E$ has numerical index $\leq \frac12.$ Taking a look at the proof of Theorem~\ref{t numerical index non-commutative JB*-triples} we observe that just a very particular case of non-commutative JB$^*$-triples remains undetermined. The following open problem covers the unique remaining case: 

\begin{problem}\label{problem non-commutative JB*-triple with atomic part Hilbert non-commutative} Let $E$ be a non-commutative JB$^*$-triple satisfying that the atomic part of $E^{**}$ is an (infinite) $\ell_{\infty}$-sum of complex Hilbert spaces regarded as type 1 Cartan factors and one of them has dimension greater than or equal to $2$. Does the inequality $n(E)\leq \frac12$ hold? 
\end{problem}

The previous problem is directly related to the next question.

\begin{problem}\label{problem all inner ideals are commutative} Let $E$ be a JB$^*$-triple. Suppose that for each $a\in E$ the inner ideal generated by $a$ is an associative JB$^*$-algebra. We know that this hypothesis does not suffices to conclude that $E$ is commutative, consider for example a Hilbert space regarded as a type 1 Cartan factor. However, can we establish a structure result describing the concrete form of such JB$^*$-triple $E$?  
\end{problem}

As we commented in the introduction, for each C$^*$-algebra $A$ (respectively, each JB$^*$-algebra $\A$) we have $n(A) = \frac12$ (respectively, $n(\A) =\frac12$) if and only if $A$ (respectively, $\A$) is non-commutative (see \cite{HuruyaPAMS1977} and \cite{iochum89commutativity}, respectively).  The arguments rely on considering the second dual of a C$^*$-algebra (respectively, a JB$^*$-algebra) and apply the corresponding Russo-Dye theorem and the abundance of unitaries in unital C$^*$- and JB$^*$-algebras. Our estimation in Theorems~\ref{t numerical index in JBW*-triples} and \ref{t numerical index non-commutative JB*-triples} is less fine in what concerns the lower bound. The general absence of unitary tripotents even in a JBW$^*$-triple discards the use of the classical arguments, and points out the interest and difficulty of the following open question which closes this note. 

\begin{problem}\label{problem bound below 1/2} Let $E$ be a JB$^*$-triple. Does the inequality $n(E)\geq \frac12$ hold?  
\end{problem}

\smallskip

\textbf{Acknowledgements}

Both authors were supported by Junta de Andalucía grants FQM375 and PY20$\underline{\ }$00255, and grant PID2021-122126NB-C31 funded by MCIN/AEI/10.13039/501100011033 and by “ERDF A way of making Europe”. First author supported by grant FPU21/00617 at University of Granada founded by Ministerio de Universidades (Spain). Second author supported by the IMAG--Mar{\'i}a de Maeztu grant CEX2020-001105-M/AEI/10.13039/501100011033.

\printbibliography[]

\end{document}